\documentclass{article}
\usepackage[utf8]{inputenc}

\usepackage{setspace}
\usepackage{amsmath}
\usepackage{amssymb}
\usepackage{fullpage}
\usepackage{graphicx}
\usepackage{imakeidx}
\usepackage{blindtext}
\usepackage{amsthm}
\usepackage{caption}
\usepackage{subcaption}
\usepackage{parskip}
\usepackage{float}
\usepackage{url}
\usepackage{geometry}
\usepackage{authblk}
\usepackage{xcolor}
\usepackage{tikz-cd}
\usepackage[colorlinks=true]{hyperref}

\definecolor{imperialBlue}{RGB}{0, 62, 116}
\definecolor{imperialBrick}{RGB}{165,25,0}
\definecolor{imperialProcess}{RGB}{0,133,202}
\definecolor{imperialGreen}{RGB}{2,137,59}
\definecolor{imperialRed}{RGB}{221,37,1}
\definecolor{imperialOrange}{RGB}{210,64,0}
\definecolor{imperialBlue2}{RGB}{0,110,175}
\definecolor{imperialTangerine}{RGB}{236,115,0}
\definecolor{imperialPurple}{RGB}{101,48,152}
\definecolor{imperialLime}{RGB}{196,214,0}
\definecolor{imperialKermit}{RGB}{102,164,10}
\hypersetup{
    colorlinks,
    citecolor=imperialGreen,
    filecolor=imperialBlue,
    linkcolor=imperialBlue,
    urlcolor=imperialProcess
}

\usepackage{booktabs} 

\newtheorem{theorem}{Theorem}[section]
\newtheorem{proposition}{Proposition}[section]

\newtheorem{definition}{Definition}[section]
\newtheorem{remark}{Remark}[section]

\newcommand{\R}{\mathbb{R}}

\newcommand{\N}{\mathbb{N}}
\newcommand{\X}{\mathcal{X}}
\newcommand{\NN}{{\mathbb{N}_0}}
\newcommand{\Rtn}{(\mathbb{R}^2)^\NN}
\newcommand{\D}{\mathcal{D}}
\newcommand{\Pcal}{\mathcal{P}}

\newcommand{\dcd}{\Delta_{c,d}}
\newcommand{\dab}{\Delta_{a,b}}

\newcommand{\intopi}{\int_0^{\pi/2}}
\newcommand{\hypf}{{}_0F_1}

\newcommand{\sigmahat}{\Hat{\sigma}}
\newcommand{\tauhat}{\Hat{\tau}}
\newcommand{\Lambdahat}{\Hat{\Lambda}}
\newcommand{\sigmatil}{\Tilde{\sigma}}
\newcommand{\tautil}{\Tilde{\tau}}
\newcommand{\ptil}{\Tilde{p}}
\newcommand{\qtil}{\Tilde{q}}

\newcommand{\phat}{\Hat{p}}
\newcommand{\qhat}{\Hat{q}}
\newcommand{\khat}{\Hat{k}}

\newcommand{\CC}{K}
\newcommand{\kxy}{k_{x,y}}

\newcommand{\sigmahatch}{\sigmahat^{ch}}

\newcommand{\tauhatch}{\tauhat^{ch}}

\newcommand{\sigmanew}{\sigma_{new}}
\newcommand{\taunew}{\tau_{new}}

\usepackage[style=alphabetic, giveninits = true, maxbibnames = 99]{biblatex}
\title{Numerical Schemes for Signature Kernels}
\author[1]{Thomas Cass}
\author[1,\dag]{Francesco Piatti}
\author[ ]{Jeffrey Pei}

\affil[1]{Department of Mathematics, Imperial College London}
\affil[$\dag$]{Corresponding author, email: francesco.piatti19@imperial.ac.uk}

\date{}

\begin{document}
\maketitle
\setstretch{1.25}

\begin{abstract}
    \noindent Signature kernels have emerged as a powerful tool within kernel methods for sequential data. In the paper ``The Signature Kernel is the solution of a Goursat PDE" \cite{salvi2021signature}, the authors identify a kernel trick that demonstrates that, for continuously differentiable paths, the signature kernel satisfies a Goursat problem for a hyperbolic partial differential equation (PDE) in two independent time variables. While finite difference methods have been explored for this PDE, they face limitations in accuracy and stability when handling highly oscillatory inputs. In this work, we introduce two advanced numerical schemes that leverage polynomial representations of boundary conditions through either approximation or interpolation techniques, and  rigorously establish the theoretical convergence of the polynomial approximation scheme. Experimental evaluations reveal that our approaches yield improvements of several orders of magnitude in mean absolute percentage error (MAPE) compared to traditional finite difference schemes, without increasing computational complexity. Furthermore, like finite difference methods, our algorithms can be GPU-parallelized to reduce computational complexity from quadratic to linear in the length of the input sequences, thereby improving scalability for high-frequency data. We have implemented these algorithms in a dedicated Python library, which is publicly available at: \url{https://github.com/FrancescoPiatti/polysigkernel}.
\end{abstract}

\section{Introduction}

Kernel methods are powerful tools in machine learning and have been widely applied in areas such as image classification and structured data analysis, achieving significant impact. They form the backbone of several machine learning techniques, including Support Vector Machines \cite{hearst1998support} and Gaussian processes \cite{rasmussen2003gaussian}, where they allow the use of high-dimensional feature sets without some of the disadvantages associated with working in high-dimensional spaces.  A kernel $\kappa:\X\times\X\to \R$ is a positive semidefinite symmetric function defined on a set $\X$. Each such kernel has associated to it a reproducing kernel Hilbert space (RKHS) $\mathcal{H}$ consisting of functions on $\X$ and containing as a dense subspace the linear span of functions $\left\{\kappa_{x}: x \in \X \right\}$ with $\kappa_{x}:=\kappa(x,\cdot):y \rightarrow \kappa (x,y)$. If $\X$ is a topological space, then a universal kernel is one for which $\mathcal{H}$ is dense in an appropriate space of continuous functions. In many applications, including those mentioned above, the kernel-based approach allows one to work with high-dimensional feature mappings without explicitly constructing those features, as only the kernel evaluations $\{\kappa(x_{i},x_{j}):i,j=1,...,n\}$ are required for a finite sample $\{x_{i}:i=1,...,n \}$.

Sequential data are inherently high-dimensional and require efficient data representation and suitable methods to mitigate the impact of dimensionality. Kernel methods tailored to ordered inputs have become increasingly widely used in fields such as finance \cite{nikan2024localized}, speech recognition \cite{may2019kernel}, and bioinformatics \cite{zou2021mk, qi2021spectral, haywood2021kernel}. Popular kernel examples include the Dynamic Time Warping (DTW) Kernel \cite{lei2007study} and the Global Alignment Kernel (GAK) \cite{cuturi2007kernel, cuturi2011fast}. An alternative approach introduced by \cite{kiraly2019kernels}, is based on the (ordinary) signature kernel, a universal kernel on the space of unparameterised paths. This uses the signature transform to represent a continuous path of finite length $\gamma$ as a hierarchy of coefficients that are determined by the collection of all iterated integrals of $\gamma$. If $\gamma$ takes values in a Hilbert space $V$, then its signature $S(\gamma)$ takes value in the free algebra of tensor series over $V$ \cite{reutenauer2003free}:
\begin{equation*}
    T((V))=\left\{\mathbf{A}=(a^0,a^1,\dots)\ |\ a^0\in\R,\forall k>0, a^k\in V^{\otimes k}\right\}.
\end{equation*}
where $\otimes$ denotes the (classical) tensor product of vector spaces. 
\begin{definition}
    Let $V$ be a Banach space, and let $\gamma:[0,t]\to V$ be a continuous path of finite $p$-variation \cite{lyons2007differential}, with $p<2$. Then the signature transform $S(\gamma)_{[0,t]}$ of the path $\gamma$ is defined as the following collection of iterated integrals:
    \begin{equation*}
        S(\gamma)_{[0,t]}=\left(1, \underset{0<u_1<t}{\int}dx_{u_1},\dots, \underset{0<u_1<\dots<u_k<t}{\int\dots\int}dx_{u_1}\otimes\dots\otimes dx_{u_k},\dots\right)
    \end{equation*}
    or, equivalently, as the unique solution to the following controlled differential equation
    \begin{equation*}
         dS(\gamma)_t=S(\gamma)_t \otimes d\gamma_t,\qquad S(\gamma)_0=1 
    \end{equation*}
\end{definition}
 It is known that the map $\gamma \mapsto S(\gamma)$ is a faithful representation of the group of unparameterised paths \cite{hambly2010uniqueness}. By equipping $T((V))$ with an inner product, we can define a kernel on the space by \[
 k_{x,y} = \langle S(x),S(y) \rangle_{T((V))}
 \]
When the inner product $\langle\cdot,\cdot\rangle_{V^{\otimes k}}$ on each $V^{\otimes k}$ is chosen as the Hilbert-Schmidt inner product induced by $\langle\cdot,\cdot\rangle_V$, the resulting inner product on $T((V))$ can be defined by linearity as $\langle\cdot,\cdot\rangle_{T((V))}=\sum_{k=0}^\infty\langle v_k,w_k\rangle_{V^{\otimes k}}$. 
This choice not only yields a universal kernel but also connects the signature kernel to a Goursat problem \cite{salvi2021signature}. While alternative inner products (e.g., weighted inner products \cite{cass2024weighted}) could be considered, they are not the focus of this paper.

The effectiveness of the signature kernel has been demonstrated in numerous studies across a variety of problem domains. Its applications include time series classification and regression \cite{salvi2021signature, gu2024transportation}, Bayesian modeling \cite{lemercier2021siggpde, toth2020bayesian}, theoretical deep learning \cite{fermanian2021framing}, and generative modeling \cite{dyer2022approximate, chevyrev2018signature, issa2024non}. Furthermore, it has been successfully applied to distribution regression \cite{lemercier2021distribution, salvi2021higher, cochrane2021sk} and numerical analysis to solve path-dependent PDEs \cite{pannier2024path}. More recently, its connections to the scaling limits of residual neural networks have been explored \cite{cirone2023neural}, further showcasing its versatility and theoretical importance.

On the other hand, the practical implementation requires a way of evaluating the inner product efficiently. The first attempt was made by \cite{kiraly2019kernels}, who proposed a truncated approximation, replacing the full signature with its truncated version up to level $N$, leveraging the factorial decay of the signature terms \cite{lyons2007differential}. Specifically, they approximated the kernel as $\langle S(\gamma), S(\sigma)\rangle\approx\langle S^N(\gamma), S^N(\sigma)\rangle$ where $S^N(\cdot)$ lies in the truncated tensor algebra $T_N(V)=\bigoplus_{k=0}^NV^{\otimes k}$. The direct evaluation of truncated signatures presents a major drawback: the number of coefficients in the signatures grows geometrically with the dimension of the time series, making it computationally expensive for large $N$. To mitigate this limitation, they introduced Horner-type schemes that exploit dynamic programming and low-rank approximations, significantly improving the computational efficiency of evaluating the truncated signature kernel. 

Another approach to approximate the truncated signature kernel is to consider Random Fourier Signature Features \cite{toth2023random}.
This method extends the framework of Random Fourier Features \cite{rahimi2007random}, which relies on Bochner’s theorem  stating that any shift-invariant kernel can be expressed as the Fourier transform of a probability measure, to sequential data via signature transforms. By leveraging this framework, the signature kernel is approximated using randomized Fourier basis functions. Although this method still requires truncation of the signature, it offers a significant computational advantage as its complexity scales linearly with the length of the input paths, making it efficient for large-scale datasets.

The untruncated signature kernel was first identified by \cite{salvi2021signature}, who demonstrated that for differentiable paths $x$ and $y$, the kernel $k_{x,y}$ satisfies a linear hyperbolic Goursat partial differential equation.

\begin{theorem}[The signature kernel PDE \cite{salvi2021signature}]\label{sigpdetheorem}
Let $I=[u,u']$ and $J=[v,v']$ be two compact intervals and let $x \in C^1(I, V)$ and $y \in C^1(J, V)$ be differentiable functions. Then, the signature kernel $\kxy$ is a solution of the following linear, second order, hyperbolic PDE:
\begin{equation}\label{sigpde}
    \frac{\partial^2 \kxy}{\partial s \partial t}=\left\langle\dot{x}_s, \dot{y}_t\right\rangle_V \kxy, \qquad \kxy(u,\cdot)=\kxy(\cdot, v)=1
\end{equation}
where $\dot{x}_s$ and  $\dot{y}_t$ are the derivatives of $x$ and $y$ at time $s$ and $t$ respectively.
\end{theorem}
The numerical scheme proposed by the authors employs finite differences methods on a refined grid (Definition \ref{dyadic scheme} in this paper), interpolating time series data to improve accuracy. However, this approach can require a high-resolution grid and tends to become numerically unstable for paths at large $p$-variation.

Rough path theory \cite{lyons1998differential} provides a natural framework to handle highly oscillatory paths, mitigating the numerical instability of finite difference methods. Additionally, it extends the computation of signature kernels to paths that may not have finite length. \cite{lemercier2024high} demonstrated that the signature kernel for smooth rough paths \cite{bellingeri2022smooth} satisfies a system of coupled PDEs. Their approach approximates input paths using piecewise log-linear paths (also known as piecewise abelian paths \cite{flint2015pathwise}), parameterized by a time partition and an approximation degree. Degree-1 approximations reduce the problem to solving a PDE similar to \eqref{sigpde}, while higher degrees result in systems of coupled PDEs with coefficients derived from log-signatures over the partition intervals. While this method enhances stability, it introduces hyperparameters and scalability challenges, particularly for large datasets.

\subsection*{Our contribution}

In this paper, we study the untruncated signature kernel and we propose two novel numerical schemes to solve PDE \eqref{sigpde} for piecewise linear input paths. First, we focus on the solution of the PDE on a general rectangle $[a,b]\times[c,d]$ where $\left\langle\dot{x}_s, \dot{y}_t\right\rangle_V=C$ is constant for $s\in[a,b]$ and $t\in[c,d]$. We show that given analytical initial boundary conditions in the form of power series expansions $k(s, c)=\sum_{n\in\NN}p_n(s-a)^n$ and $k(a,t)=\sum_{n\in\NN}q_n(t-c)^n$  the coefficients of the power series that describe the solution on the opposite edges of the rectangle, i.e. $k(s, d)=\sum_{n\in\NN}\ptil_n(s-a)^n$ and $k(b,t)=\sum_{n\in\NN}\qtil_n(t-c)^n$ can be derived analytically 
\begin{align*}
    \ptil_{n} &= \sum_{k=0}^{n} p_k \frac{(C (t-c))^{n-k}\ k!}{(n-k)!\ n!} + \sum_{k=1}^{\infty} q_k \frac{ C^n(t-c)^{n+k}\ k!}{(n+k)!\ n!} \\
    \qtil_{n} &= \sum_{k=0}^{n} q_k \frac{(C (s-a))^{n-k}\ k!}{(n-k)!\ n!} + \sum_{k=1}^{\infty} p_k \frac{ C^n(s-a)^{n+k}\ k!}{(n+k)!\ n!}
\end{align*}
This result, formalized and proven in Proposition \ref{approxsolution}, provides an analytical tool to propagate the boundary conditions across the grid and compute the signature kernel recursively. Building on this result, we develop two numerical schemes that leverage polynomials to approximate \eqref{sigpde}. The first scheme, referred to as the polynomial approximation scheme, truncates the power series above at order $N$. For this scheme, we derive an explicit error bound of the form:
\begin{equation*}
        |k(s_{\ell_x},t_{\ell_y})- \khat^N(s_{\ell_x},t_{\ell_y})|\le  f(\ell_x+\ell_y - 1)\frac{(\CC\Delta^2)^{N+1}}{[(N+1)!]^2}\left[\frac{(\ell_x+\ell_y - 1)^{N+2}}{N+2}+(\ell_x+\ell_y - 1)^{N+1}\right]
\end{equation*}
        where $\Delta=\max_{i,j}\{s_{i+1}-s_i,t_{j+1}-t_j\}$, $ \sup\left|\left\langle\dot{x}_s, \dot{y}_t\right\rangle\right|\le \CC$,  $f:\N\to\R$ is a function defined as $f(k):=\prod_{m=0}^kI_0(2\sqrt{m\CC}\Delta)$, and $\ell_x$, $\ell_y$ are the number of linear segments of $x$ and $y$.

The second scheme, termed the polynomial interpolation scheme, approximates the solution through polynomial interpolation and achieves faster convergence, particularly for oscillatory paths. By leveraging these polynomial-based approaches, we significantly reduce the need for high-resolution grids required by finite differences methods, improving numerical stability, especially for paths with complex oscillatory behavior.

Both numerical schemes have been integrated into an open-source software library, which also includes tools for practical applications of the signature kernel, such as Maximum Mean Discrepancy (MMD) computation and statistical hypothesis testing.

\subsection*{Outline}
Section \ref{sec:pre} provides an overview of the finite differences scheme introduced by \cite{salvi2021signature} for solving the signature kernel PDE \eqref{sigpde}, which serves as a baseline for comparison with our methods. Section \ref{sec:main} is divided into three parts. In Section \ref{sec:31}, we solve the PDE \eqref{sigpde} on a grid recursively by explicitly expressing the solution as a power series and we establish bounds on the coefficients. Building on these results, we introduce two polynomial-based numerical schemes: the polynomial approximation scheme (Section \ref{sec:polyapprox}) and the polynomial interpolation scheme (Section \ref{sec:polyinterp}). Finally, in Section \ref{sec:exp}, we assess the performance of our schemes in terms of accuracy and computational efficiency. 

\section{Finite Differences Scheme}\label{sec:pre}

In this section, we briefly review the finite differences scheme proposed by \cite{salvi2021signature} for solving the signature kernel PDE \eqref{sigpde} when the input paths $x$ and $y$ are piecewise linear. Throughout this paper, we use the following notation: \(C(I, \R)\) denotes the class of continuous functions defined on a compact interval $I$, \(C^k(I, \R)\) represents the class of $k$-times differentiable functions, and \(C^\infty(I, \R)\) denotes the class of smooth functions.

Suppose $x:I\to V$ and $y:J\to V$ are piecewise linear paths on the partitions $\D_I = \{u_{0} < \dots < u_{m}\}$ and $\D_{J} = \{v_{0} < \dots < v_{n}\}$ respectively. To simplify notation also assume that $V=\R^d$. Then PDE \eqref{sigpde} becomes
\begin{equation}\label{PDERect1}
    \frac{\partial^2 \kxy}{\partial s \partial t}=C \kxy
\end{equation}
on each domain $\D_{ij}=\{(s,t) | u_i \le s \le u_{i+1},v_j \le t \le v_{j+1}\}$ where $C=\left\langle\dot{x}_s, \dot{y}_t\right\rangle_V=\,\textrm{const.}$.

In \cite{salvi2021signature}, explicit and implicit finite differences schemes are derived to approximate the solution of PDE \eqref{PDERect1}, which can be expressed in integral form as:
\begin{equation}\label{integralpde}
    \kxy(s,t)=\kxy(s,v) + \kxy(u, t) - \kxy(u,v) + C\int_u^s\int_v^t \kxy(r,w)dwdr
\end{equation}
for $(s,t),(u,v)\in\D_{ij}$, with $u\le s$ and $v\le t$. By approximating the double integral in \eqref{integralpde} using the trapezoidal rule, the following finite differences scheme is obtained:
\begin{align*}
    \kxy(s,t)&\approx\kxy(s,v) + \kxy(u, t) - \kxy(u,v) + \frac{1}{2}C\big(\kxy(s,v) + \kxy(u, t)\big)(s-u)(t-v)
\end{align*}

\begin{remark}\label{higherfd}
    Higher-order quadrature methods can provide more accurate approximations for the double integral in \eqref{integralpde}. By including all four values of the kernel, they derive an implicit scheme:
    \begin{align*}
        \kxy(s,t)&\approx\kxy(s,v) + \kxy(u, t) - \kxy(u,v)  \\ &+\frac{1}{2}C\big(\kxy(u,v) + \kxy(s,v) + \kxy(u, t) + \kxy(s,t)\big)(s-u)(t-v) \nonumber
    \end{align*}
    They also introduce a second-order explicit scheme:
    \begin{align}\label{secondexplicitscheme}
        \kxy(s,t)\approx\kxy(s,v) &+ \kxy(u, t) - \kxy(u,v) + \frac{1}{2}C\big(\kxy(s,v) + \kxy(u, t)\big)(s-u)(t-v) \\
        &+\frac{1}{12}C^2\big(\kxy(s,v) + \kxy(u, t) + \kxy(u,v)\big)(s-u)^2(t-v)^2 \nonumber
    \end{align}
\end{remark}

Building on these foundations, a dyadic refinement scheme is proposed, which subdivides each rectangle of the initial grid $\D_I\times\D_J$ into smaller rectangles via linear interpolation of the time series. This refinement results in a new grid $P_\lambda$, derived from the $\lambda$-order dyadic refinement of the initial partition $\Pcal_0 = \D_I \times \D_J$. 

\begin{definition}[Dyadic finite differences approximation scheme, \cite{salvi2021signature}]
    \label{dyadic scheme}
     Let $\Pcal_{0} = \D_{I} \times \D_{J}$ and for integers $\lambda \ge 1$ define $\Pcal_{\lambda}$ to be the $\lambda$-order dyadic refinement of $\Pcal_{0}$, that is,
    \begin{align*}
        \Pcal_{\lambda}
        \cap
        \left(
            [u_{i}, u_{i+1}]\times
            [v_{j}, v_{j+1}]
        \right)
        &=
        \left\{
            \left(
                u_{i}
                +
                k2^{-\lambda}(u_{i+1} - u_{i}),
                v_{j}
                +
                l2^{-\lambda}(v_{j+1} - v_{j})
            \right)
            \;:\;
            0 \leq k,l \leq 2^{\lambda}
        \right\}
    \end{align*}
    On the grid $\Pcal_{\lambda} = \{(s_{i}, t_{j}) : 0\le i\le 2^{\lambda} m, \ 0 \le j \le 2^{\lambda}n\}$, the first-order finite differences scheme $\khat = \khat^{\lambda}$ is defined as
    \begin{align*}
        \khat(s_{i+1},t_{j+1})
        =
        \begin{aligned}[t]
            &
            \khat(s_{i+1},t_{j})
            +
            \khat(s_{i},t_{j+1})
            -
            \khat(s_{i},t_{j})
            \\
            &+
            \frac{1}{2}
            \left\langle
                x_{s_{i+1}} - x_{s_{i}},
                y_{t_{j+1}} - y_{t_{j}}
            \right\rangle
            \left(
                \khat(t_{i+1},v_{j})
                +
                \khat(t_{i},v_{j+1})
            \right).
        \end{aligned}
    \end{align*}
    where $\khat(s_{0}, \cdot) = \khat(\cdot, t_{0}) = 1$. Similarly, one can define higher-order implicit and explicit schemes (Remark \ref{higherfd}) on the same $\Pcal_\lambda$.
\end{definition}

As the dyadic order $\lambda$ increases, the numerical solution converges to the true solution, yet significantly increasing the computational complexity. Theorem 3.5 of \cite{salvi2021signature} proves that the error of the finite difference approximation scheme decreases at a rate proportional to $2^{-2\lambda}$.

\section{Polynomial-based schemes for the Goursat PDE}\label{sec:main}

In this section, first we derive an analytical solution for the signature kernel of two piecewise-linear paths in the form of a power series with coefficients that can be computed recursively, and  we establish bounds for these coefficients. We then introduce two numerical schemes which leverage polynomials to approximate the  signature kernel: the first scheme approximates the signature kernel by truncating the aforementioned power series at order $N$, while the second scheme employs polynomial interpolation to numerically solve the Goursat PDE.

\subsection{Recursive solution for the Goursat PDE}\label{sec:31}

Consider the following case of the Goursat problem \eqref{sigpde} on the rectangular domain $\mathcal{D}_C=\left\{(s, t) \mid a \leq s \leq b, c \leq t \leq d\right\}$:
\begin{equation}\label{goursatonrectangle}
    \frac{\partial^2 k}{\partial s \partial t}=C k
\end{equation}
where $C\in\R/\{0\}$ is constant.
\begin{theorem}[\cite{polyanin2001handbook}]\label{besselsolution}
     The solution $k$ to the Goursat problem \eqref{goursatonrectangle} on the domain $\D_C$  with differentiable boundary data $k(s, c)=\sigma(s) \in C^1([a,b],\R^d)$ and $k(a, t)=\tau(t)\in C^1([c,d],\R^d)$ can be expressed as
\begin{equation}\label{eq:integralsolution}
    k(s, t)=k(a, c) R(s-a, t-c)+\int_a^s \sigma^{\prime}(r) R(s-r, t-c) d r+\int_c^t \tau^{\prime}(r) R(s-a, t-r) d r
\end{equation} for $(s, t) \in \mathcal{D}_C$. Here, $R$ is the the Riemann function defined as $R(w, z):=J_0\left(2 i \sqrt{C w z}\right)$ for $w,z \geq 0$, with $J_0$ denoting the zero-order Bessel function of the first kind. 
\end{theorem} 
Equation \eqref{eq:integralsolution} can also be written in terms of the zero-order modified Bessel function of the first kind.
\begin{align*}
k(s, t)&=\sigma(a) I_0\left(2 \sqrt{C(s-a)(t-c)} \right) +\int_a^s \sigma^{\prime}(r) I_0\left(2 \sqrt{C(s-r)(t-c)}\right) d r +\int_c^t \tau^{\prime} (r) I_0 \left(2 \sqrt{C(s-a)(t-r)}\right) d r 
\end{align*}
\begin{remark}\label{mbf}
The $n$-th order modified Bessel function is defined as $I_n(z):=i^{-n} J_n(i z)$.
It directly follows from the series expansion of $J_n(2 z)$ \cite{andrews1999special} that
\begin{equation*}
    I_n(2 z)=\sum_{k=0}^{\infty} \frac{z^{2 k+n}}{k!(n+k)!}
\end{equation*}
\end{remark}

In the following sections, we leverage Theorem \ref{besselsolution} to study the signature kernel PDE \eqref{sigpde}. We start by introducing the \textit{one-step estimate}, where we solve the PDE on a single rectangle, as formulated in \eqref{goursatonrectangle}, obtaining the solution as a power series and establishing bounds on its coefficients. We then extend this to the \textit{multi-step estimate}, deriving a recursive formulation for the exact signature kernel of piecewise linear paths and proving bounds on the coefficients that characterize this recursive solution.

\subsubsection*{One-step estimate}

Consider the Goursat problem on the rectangle $[a,b]\times[c,d]$ as in \eqref{goursatonrectangle}. Assume that the initial boundary conditions $\sigma$ and $\tau$ are (real) analytic functions at $s=a$, $t=c$ respectively, and are determined by their power series
\begin{equation}\label{initialps}
    \sigma(s) = k(s, c) = \sum_{n\in\NN} p_n (s-a)^n \qquad \tau(t) = k(a, t) = \sum_{n\in\NN} q_n (t-c)^n
\end{equation}
with the constrain $\sigma(a)=\tau(c)$, as both represent the value $k(a,c)$. Furthermore, assume that these power series converge uniformly on $[a,b]$ and $[c,d]$ respectively.

We are interested in understanding the functions:
\begin{equation}\label{finalps}
        \sigmatil(s) = k(s,d) = \sum_{n\in\NN} \ptil_n (s-a)^n \qquad \tautil(t) = k(b, t) =  \sum_{n\in\NN} \qtil_n (t-c)^n
\end{equation}
which represent the solution of the PDE on the top and right edges of $[a,b]\times[c,d]$ respectively.

Specifically, we want:
\begin{enumerate}
    \item to determine how $(\ptil, \qtil):=(\ptil_n,\qtil_n)_{n\in\NN}\in\Rtn$ is derived from $(p,q)$, i.e. what is the function 
    \begin{equation}\label{lambda func}
        \Lambda_{C;\,[a,b]\times[c,d]} \,:\, (\R^2)^{\NN} \to  (\R^2)^{\NN}\qquad \textrm{such that}\qquad(\ptil, \qtil) =\Lambda_{C,\, [a,b]\times[c,d]}(p,q)
    \end{equation}
    \item to understand how we may compare the magnitude of the terms in the sequence $(\ptil,\qtil)$ with those in $(p,q)$. To this end we define a norm on $\Rtn$ by
    \begin{equation}\label{norm}
        \|(p^1,p^2)\|_{\gamma} = \max_{i=1,2}\left\{\max_{n\in\NN}\left\{\left|\frac{(n!)^2p^i_n}{(\gamma |C|\Delta)^n}\right|\right\}\right\}
    \end{equation}
    where $\Delta=\max\{b-a,d-c\}$ and $1\le\gamma\in \N$.
\end{enumerate}

\begin{proposition}\label{approxsolution}
    Consider the Goursat PDE on the rectangle $[a,b]\times[c,d]$ as in \eqref{goursatonrectangle}, with initial boundaries condition given by \eqref{initialps}, then $\forall n\in\NN$ the coefficients of \eqref{finalps} are given by
    \begin{equation}\label{powerseries1}
        \ptil_{n} = \sum_{k=0}^{n} p_k \frac{(C \dcd)^{n-k}\ k!}{(n-k)!\ n!} + \sum_{k=1}^{\infty} q_k \frac{ C^n\dcd^{n+k}\ k!}{(n+k)!\ n!}
    \end{equation}
    and, symmetrically
    \begin{equation}\label{powerseries2}
        \qtil_{n} = \sum_{k=0}^{n} q_k \frac{(C \dab)^{n-k}\ k!}{(n-k)!\ n!} + \sum_{k=1}^{\infty} p_k \frac{ C^n\dab^{n+k}\ k!}{(n+k)!\ n!}
    \end{equation}
    where $\dab=(b-a)$ and $\dcd=(d-c)$.
\end{proposition}
Refer to Appendix \ref{app:A1} for the proof of the above proposition. 

Proposition \ref{approxsolution} allows us to explicitly define the linear map $\Lambda_{C;\,[a,b]\times[c,d]}: (\R^2)^{\NN}\to(\R^2)^{\NN}$ introduced in \eqref{lambda func}. As we are interested in the behavior of the resulting coefficients, it is natural to consider $\Lambda_{C;\,[a,b]\times[c,d]} \,:\, ((\R^2)^{\NN}, \|\cdot\|_\gamma) \to  ((\R^2)^{\NN}, \|\cdot\|_{\gamma+1})$. The following proposition shows that this map is bounded

\begin{proposition}\label{prop:onestepestimate}
    
    Consider the linear map $\Lambda_{C;\,[a,b]\times[c,d]} \,:\, ((\R^2)^{\NN}, \|\cdot\|_\gamma) \to  ((\R^2)^{\NN}, \|\cdot\|_{\gamma+1})$, as in \eqref{lambda func}, for some $1\le\gamma\in\N$, where $\|\cdot\|_\gamma$ is defined in \eqref{norm}. Then $\Lambda_{C;\,[a,b]\times[c,d]}$ is bounded and 
    \begin{equation*}
        \|\Lambda_{C;\,[a,b]\times[c,d]}\|_{op}\le I_0(2\sqrt{\gamma |C|}\Delta)
    \end{equation*}
    where $\|\cdot\|_{op}$ is the operator norm and $\Delta=\max\{b-a,d-c\}$. 
\end{proposition}

\paragraph{Proof.} Fix $1\le\gamma\in\N$ and assume $\|(p,q)\|_\gamma\in\R$. Then observe that:
\begin{align*}
    |(n!)^2\ptil_n|&\le \sum_{k=0}^{n}\binom{n}{k} (k!)^2p_k (|C|\dcd)^{n-k} + n!(|C|\Delta_{c,d})^n\sum_{k=1}^\infty \frac{(k!)^2q_k\Delta_{c,d}^k}{(n+k)!k!} \\
    &\le(|C|\Delta)^{n}\|(p,q)\|_\gamma\sum_{k=0}^n\binom{n}{k}\gamma^k + (|C|\Delta)^n\|(p,q)\|_\gamma\sum_{k=1}^\infty \frac{(|C|\gamma)^k\Delta^{2k}}{(k!)^2} \\
    &\le (|C|\Delta)^n\|(p,q)\|_\gamma(\gamma+1)^n + (|C|\Delta)^n\|(p,q)\|_\gamma(I_0(2\sqrt{\gamma |C|}\Delta)-1)
\end{align*}
where in the second inequality we have used $(n+k)!\ge n!k!$.
This implies that: 
\begin{align*}
    \left|\frac{(n!)^2\ptil_n}{(C\Delta)^n(\gamma+1)^n}\right| & \le \|(p,q)\|_\gamma + \|(p,q)\|_\gamma \frac{(I_0(2\sqrt{\gamma |C|}\Delta)-1)}{(\gamma+1)^n} \\
    & \le \|(p,q)\|_\gamma I_0(2\sqrt{\gamma |C|}\Delta) 
\end{align*} 
A similar argument holds for $\qtil_n$, leading to $\|(\ptil,\qtil)\|_{\gamma+1}\le\|(p,q)\|_\gamma I_0(2\sqrt{\gamma |C|}\Delta)$.\qed

\begin{remark}
    The bound on \(\|(\ptil, \qtil)\|_{\gamma+1}\) guarantees that the power series expansions in \eqref{finalps} converge uniformly on the intervals $[a, b]$ and $[c, d]$, respectively, which is essential for a recursive application of Proposition \eqref{prop:onestepestimate}.
\end{remark}

\subsubsection*{Multi-step estimate}

In the one-step estimate we derived analytical formulas to solve PDE \eqref{goursatonrectangle} on a rectangle and we derived a bound for $\|(\ptil, \qtil)\|_{\gamma+1}$ in terms of $\|(p,q)\|_{\gamma}$. We now extend this result by developing a recursive method to obtain a similar bound for the coefficients of the power series describing the solution of PDE \eqref{sigpde} over a grid induced by piecewise linear paths. Notably, Theorem \ref{sigpdetheorem} remains well-defined even when the input paths are only piecewise $C^1$.

Suppose we are given the following data:
\begin{itemize}
    \item Two partitions $\D_{I}=\{s_i : i=0, 1, \dots, l_x\}=\{0=s_0<s_1<\dots<s_i<\dots<s_{\ell_x}\}$ and, similarly, $\D_{J}=\{t_j : j=0, 1, \dots, \ell_y\}$. 
    \item A subset $\{C_{ij} \in \R : i=0,1,\dots,\ell_x-1, \ j=0,1,\dots,\ell_y-1\}$
\end{itemize}
We aim to solve 
\begin{equation*}
    \frac{\partial^2k}{\partial s\partial t} = C(s,t)k \qquad \textrm{on}\quad [0,s_{\ell_x}]\times[0,t_{\ell_y}]
\end{equation*}
subject to initial boundary conditions $k(0,\cdot)= k(\cdot, 0)= 1$, where $C(s,t)\equiv C_{ij}$ for $(s,t) \in (s_i, s_{i+1}]\times (t_j, t_{j+1}]$.

To extend the Proposition \ref{prop:onestepestimate} to a grid-based approach, we introduce the sequences $p^{ij}=(p^{ij}_n)_{n\in\NN}$ and $q^{ij}=(q^{ij}_n)_{n\in\NN}$ as the coefficients of the power series defining the functions $k(\cdot, t_j):[s_i,s_{i+1}]\to\R$ and $k(s_i, \cdot):[t_i,t_{i+1}]\to\R$ respectively, i.e.
\begin{align}
    k(s,t_j) &= \sum_{n\in\NN}p_n^{ij}(s-s_i)^n \qquad s \in [s_i,s_{i+1}] \label{ps1} \\
    k(s_i,t) &= \sum_{n\in\NN}q_n^{ij}(t-t_j)^n \qquad t \in [t_j,t_{j+1}] \label{ps2} 
\end{align}

To improve clarity, we introduce additional notation to represent the power series coefficients associated with the boundary solutions of a each rectangle. Furthermore, we extend the operator $\Lambda_{C;\,[a,b]\times[c,d]}$, defined in \eqref{lambda func}, to each rectangle in the grid $\D_I\times\D_J$.

\begin{definition}\label{lambda}
    The sequences of coefficients $(p^{ij}_n)_{n\in\NN}$ and $(q^{ij}_n)_{n\in\NN}$, which characterise power series \eqref{ps1} and \eqref{ps2} respectively, are defined recursively by 
    \begin{equation*}
        p^{i0} = q^{0j} = (1,0,0,\dots), \qquad\textrm{for}\quad i=0,1,\dots,\ell_x-1,\ j=0,1,\dots,\ell_y-1
    \end{equation*} and
    \begin{equation*}
        (p^{i,j+1}, q^{i+1,j}) := \Lambda_{C_{ij};[s_i,s_{i+1}]\times[t_j,t_{j+1}]} (p^{ij},q^{ij})
    \end{equation*}
    where $\Lambda_{C_{ij};[s_i,s_{i+1}]\times[t_j,t_{j+1}]}$ is the operator described by Proposition \ref{approxsolution}, extended to a general rectangle $[s_i,s_{i+1}]\times[t_j,t_{j+1}]$ with coefficient $C=C_{ij}$. Hence, this operator takes the coefficients $(p^{ij},q^{ij})$ of the power series describing the initial boundaries of the rectangle $[s_i,s_{i+1}]\times[t_j,t_{j+1}]$ and produces the coefficients of the power series describing the top and right boundaries, effectively solving the Goursat PDE on that rectangle.
    
    Additionally we define $\Lambda_{ij}\in\Rtn$ as
    \begin{equation*}
        \Lambda_{ij}: = (p^{i,j+1}, q^{i+1,j})\in \Rtn, \qquad i=0,1,\dots,\ell_x-1,\ j=0,1,\dots,\ell_y-1
    \end{equation*}
    where each $\Lambda_{ij}$ labels the sequence of coefficients associated with the top and right boundaries of the corresponding rectangle in the grid.
\end{definition} 

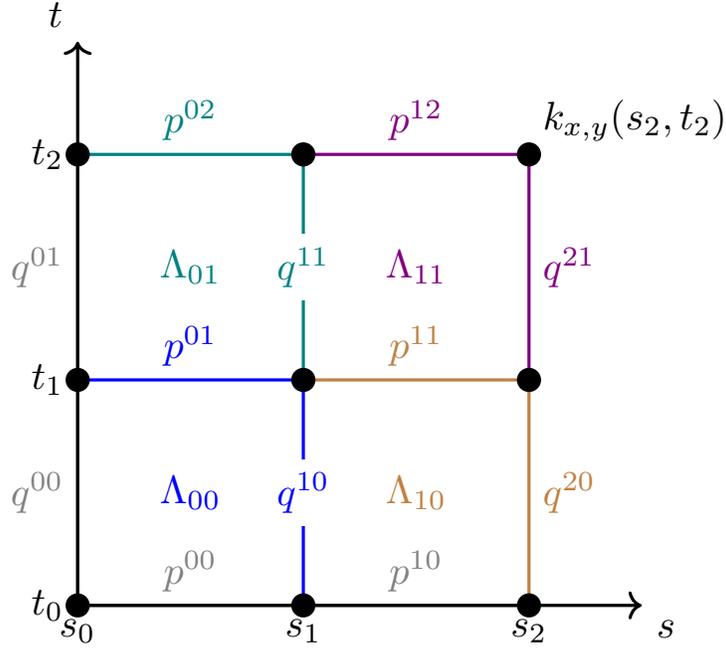
\begin{figure}[h!]
    \centering
    \scalebox{1.5}{
        \begin{tikzpicture}
            \draw[thick,gray] (0,0) -- (2,0);   
            \draw[thick,blue] (2,0) -- (2,2);   
            \draw[thick,blue] (0,2) -- (2,2);   
            \draw[thick,gray] (0,0) -- (0,2);   
            \draw[thick,gray] (2,0) -- (4,0);  
            \draw[thick,brown] (4,0) -- (4,2);  
            \draw[thick,brown] (2,2) -- (4,2);   
            \draw[thick,teal] (2,2) -- (2,4);   
            \draw[thick,teal] (0,4) -- (2,4);   
            \draw[thick,gray] (0,2) -- (0,4);   
            \draw[thick,violet] (4,2) -- (4,4); 
            \draw[thick,violet] (2,4) -- (4,4); 
            \draw[thick,->] (0,0) -- (5,0) node[anchor=north west] {$s$};
            \draw[thick,->] (0,0) -- (0,5) node[anchor=south east] {$t$};
            \draw (0,0) node[anchor=north] {$s_0$};
            \draw (2,0) node[anchor=north] {$s_1$};
            \draw (4,0) node[anchor=north] {$s_2$};
            \draw (0,0) node[anchor=east] {$t_0$};
            \draw (0,2) node[anchor=east] {$t_1$};
            \draw (0,4) node[anchor=east] {$t_2$};
            \draw (1,0) node[anchor=south] {\color{gray}$p^{00}$};
            \draw (3,0) node[anchor=south] {\color{gray}$p^{10}$};
            \draw (1,2) node[anchor=south] {\color{blue}$p^{01}$};
            \draw (3,2) node[anchor=south] {\color{brown}$p^{11}$};
            \draw (1,4) node[anchor=south] {\color{teal}$p^{02}$};
            \draw (3,4) node[anchor=south] {\color{violet}$p^{12}$};
            \draw (0, 1) node[anchor=east] {\color{gray}$q^{00}$};
            \draw (0, 3) node[anchor=east] {\color{gray}$q^{01}$};
            \draw (4, 1) node[anchor=west] {\color{brown}$q^{20}$};
            \draw (4, 3) node[anchor=west] {\color{violet}$q^{21}$};
            \draw (1,1) node[anchor= center] {\color{blue}$\Lambda_{00}$};
            \draw (3,1) node[anchor= center] {\color{brown}$\Lambda_{10}$};
            \draw (1,3) node[anchor= center] {\color{teal}$\Lambda_{01}$};
            \draw (3,3) node[anchor= center] {\color{violet}$\Lambda_{11}$};
            \node[fill=white] at (2,1) {\color{blue}$q^{10}$};  
            \node[fill=white] at (2,3) {\color{teal}$q^{11}$};  

            \filldraw[black] (0,0) circle(.1);
            \filldraw[black] (0,2) circle(.1);
            \filldraw[black] (0,4) circle(.1);
            \filldraw[black] (2,0) circle(.1);
            \filldraw[black] (2,2) circle(.1);
            \filldraw[black] (2,4) circle(.1);
            \filldraw[black] (4,0) circle(.1);
            \filldraw[black] (4,2) circle(.1);
            \filldraw[black] (4,4) circle(.1) node[anchor=south west] {$k_{x,y}(s_2, t_2)$};
        \end{tikzpicture}
    }
    \caption{Each $\Lambda_{ij}\in(\R^2)^\NN$, indexed by the grid’s rectangles, represents the power series coefficients corresponding to the top and right boundaries of that rectangle. These boundaries are colour-coded in the figure for clarity, e.g. $\Lambda_{00}:=(p^{01}, q^{10})$ and $\Lambda_{01}:=(p^{02}, q^{11})$.}
    \label{fig:grid} 
\end{figure}

Figure \ref{fig:grid} (below) illustrates the structure of the grid and the objects we have defined. 

\begin{remark}\label{symmetry1}
    The function $k(s_{\ell_x},t_{\ell_y})$, which determines the signature kernel, can be obtained evaluating power series with coefficients given by $\Lambda_{\ell_x-1,\ell_y-1}$. Further, by the symmetry of our problem, it holds that 
    \begin{equation*}
        k(s_{\ell_x}, t_{\ell_y})= \sum_{n=0}^\infty p^{\ell_x-1,\ell_y}_n(s_{\ell_x}-s_{\ell_x-1})^n = \sum_{n=0}^\infty q^{\ell_x,\ell_y-1}_n(t_{\ell_y}-t_{\ell_y-1})^n
    \end{equation*}
\end{remark}

Our goal is now to establish bounds on $\|\Lambda_{ij}\|_{\gamma}$ for all $i,j$ for which it is determined, with $1\le\gamma\in\N$. Let $\CC=\max_{i,j}|C_{ij}|$ and $\Delta=\max_{i,j}\{s_{i+1}-s_i,t_{j+1}-t_j\}$ and modify the norm \eqref{norm} such that:
\begin{equation}\label{norm2}
    \|(p^1,p^2)\|_{\gamma} = \max_{i=1,2}\left\{\max_{n\in\NN}\left\{\left|\frac{(n!)^2p^i_n}{(\gamma \CC\Delta)^n}\right|\right\}\right\}
\end{equation}
\begin{proposition}\label{multistepestimate}
    Let $\Lambda_{ij}$ be as introduced in Definition \ref{lambda}, and define the norm $\|\cdot\|_\gamma$ as in \eqref{norm2}. Then, for all $i=0,1,\dots,\ell_x-1, \ j=0,1,\dots,\ell_y-1$ it holds that 
    \begin{equation*}
        \|\Lambda_{ij}\|_{i+j+1}\le f(i+j)
    \end{equation*}
    where $f:\N\to\R$ s.t.
    \begin{equation}\label{fdefinition}
        f(k):=\prod_{m=0}^kI_0(2\sqrt{m\CC}\Delta)
    \end{equation}
    with $\CC=\max_{i,j}|C_{ij}|$ and $\Delta=\max_{i,j}\{s_{i+1}-s_i,t_{j+1}-t_j\}$.
\end{proposition}
\paragraph{Proof.}
To begin, notice that 
\begin{align*}
    p^{01}_n = \frac{C_{00}^n t_1^n}{(n!)^2},\qquad q^{10}_n = \frac{C_{00}^n s_1^n}{(n!)^2}, \qquad \forall n\in\NN
\end{align*}
which implies $\|\Lambda_{00} \|_1 = \|(p^{01}, g^{10})\|_1 \le 1$, as $C_{00}\le \CC$. By a simple induction, $\|\Lambda_{0j}\|_{j+1} \le f(j)$ and $\|\Lambda_{i0}\|_{i+1} \le f(i)$.  In fact, for $\|\Lambda_{0j}\|_{j+1}=\|(p^{0,j+1},q^{1,j})\|_{j+1}$ notice that 
\begin{equation*}
    p^{0,j+1}_n = \sum_{k=0}^{n} p^{0j}_k \frac{(C_{0j}  (t_{j+1}-t_j))^{n-k}\ k!}{(n-k)!\ n!} \qquad\textrm{and}\qquad  q^{1,j}_n=\sum_{k=0}^{\infty} p^{0j}_k \frac{ C_{0j}^n s_1^{n+k}\ k!}{(n+k)!\ n!}
\end{equation*}
where we have used $p^{0j}_0=q^{0j}_0$ and that $q^{0j}_{n}=0$ for $n\ge1$. By Proposition \ref{prop:onestepestimate}, $\|\Lambda_{0,j+1}\|_{j+2}\le\|\Lambda_{0j}\|_{j+1}I_0(2\sqrt{(j+1)\CC}\Delta)$, leading to $\|\Lambda_{0,j+1}\|_{j+2}\le f(j+1)$, by definition of $f$. Similarly we can show that, $\|\Lambda_{i0}\|_{i+1}\le f(i)$.

Next, proceed by induction on $m=i+j$. Assume $\|\Lambda_{ij}\|_{i+j+1}\le f(i+j)$ for all $i,\ j\in\N$ such that $i+j=m$ and $\Lambda_{ij}$ is defined. This implies $\|\Lambda_{ij}\|_{m+1}\le f(m)$. Exploiting the symmetry of the problem, it suffices to show  $\|\Lambda_{i+1,j}\|_{i+j+2}\le f(i+j+1)$, i.e. $\|\Lambda_{i+1,j}\|_{m+2}\le f(m+1)$.
By definition 
\begin{equation*}
    \Lambda_{i+1,j} = (p^{i+1,j+1}, q^{i+2,j}) = \Lambda_{C_{i+1,j};[s_{i+1},s_{i+2}]\times[t_j,t_{j+1}]} (p^{i+1,j},q^{i+1,j})
\end{equation*}
where $(p^{i+1,j}, q^{i+2,j-1})=\Lambda_{i+1,j-1}$ and $(p^{i,j+1}, q^{i+1,j})=\Lambda_{i,j}$.

By the inductive hypothesis $\|\Lambda_{i+1,j-1}\|_{m+1}\le f(m)$ and $\|\Lambda_{i,j}\|_{m+1}\le f(m)$, implying $\|(p^{i+1,j},q^{i+1,j})\|_{m+1}\le f(m)$. Then by Proposition \ref{prop:onestepestimate}, we obtain $\|\Lambda_{i+1,j}\|_{m+2}\le I_0(2\sqrt{(m+1)\CC}\Delta)f(m)=f(m+1)$.\qed

This result provides bounds on the coefficients of the power series describing each boundary of the grid, which in turn allows us to prove the error bounds of the numerical scheme proposed in the next section.

\subsection{Polynomial approximation numerical scheme}\label{sec:polyapprox}

We present an approximation scheme that leverages the results of Proposition \ref{approxsolution} to numerically solve the Goursat PDE \eqref{sigpde}. This scheme approximates each boundary of the grid using an order-$N$ polynomial obtained by truncating the power series \eqref{powerseries1}, \eqref{powerseries2}. This approach provides an efficient and accurate method for computing the signature kernel.

To proceed, we need to introduce additional notation to represent truncated sequences of coefficients, as well as sequences starting from an index greater than 0.
\begin{definition}
    Let $m\in\N$ and $g=(g_n)_{n\in\NN}\in\R^{\NN}$ and $h=(h_n)_{n\in\NN}\in\R^{\NN}$, then define
    \begin{equation*}
        g^{\le m}:= (g_0, g_1, \dots, g_m, 0, \dots)
    \end{equation*}
    and $(g,h)^{\le m}:=(g^{\le m}, h^{\le m})$. Similarly, define the sequence $g^{>m}$ starting from index $m+1$ as
    \begin{equation*}
        g^{>m} := g - g^{\le m}=(0, \dots, 0, g_{m+1}, g_{m+2}, \dots)
    \end{equation*}
In an obvious way we also write $g^{\ge m} = g^{>m-1}$. 
\end{definition}

Next, we define the space 
\begin{equation}\label{espace}
    E_N=\{e\in\R^\NN: e=e^{\le N}\}\equiv\{(e_0,e_1,\dots)\in \R^\NN: e_n=0\ \textrm{for}\ n>N\}
\end{equation} which encodes the order-$N$ polynomial coefficients, embedded in the space of sequences. This framework enables a more efficient comparison between the polynomial coefficients of the numerical scheme and the power series coefficients of the exact solution. Finally, we define the projection map that allows us to move from the untruncated space $\R^\NN$ to $E_N$.

\begin{definition}\label{truncmap}
     Let $m\in\N$ and $(g,h)\in(\R^2)^{\NN}$. We define the truncation operator $\pi_m : (\R^2)^\NN\to E_m\times E_m = E_m^2$ such that $\pi_m(g,h)=(g,h)^{\le m}$, where $E_m$ is defined in \eqref{espace}.
\end{definition}

With this machinery in place, we can now formally define the polynomial approximation scheme.

\begin{definition}[Polynomial approximation numerical scheme]\label{def:polyapproxscheme}
Let $\D_I =\{s_0<s_1<\dots<s_{\ell_x}\}$ and $\D_J =\{t_0<t_1<\dots<t_{\ell_y}\}$ be partitions of the compact intervals $I$ and $J$ and let $x:I\to\R^d$ and $y:J\to\R^d$ be two piecewise linear paths on $\D_I$ and $\D_J$ respectively. 
We define $\khat^N$, the solution of the order-$N$ polynomial approximation scheme to the Goursat PDE \eqref{sigpde} on the grid $\D_I\times\D_J$, as follows: 

For fixed $N\ge2$, let $\phat^{ij},\ \qhat^{ij}\in E_N$, where $E_N=\{e\in\R^\NN: e=e^{\le N}\}$, represent the coefficients of the polynomials approximating the grid edges, satisfying: 
\begin{align*}
    \khat^N(s,t_j) &= \sum_{n=0}^N \phat^{ij}_n(s-s_i)^n,\qquad s\in[s_i,s_{i+1}]\\
    \khat^N(s_i,t) &= \sum_{n=0}^N \qhat^{ij}_n(t-t_j)^n,\qquad t\in[t_j,t_{j+1}]
\end{align*}
Let  
\begin{equation*}
    C_{ij} = \frac{\langle x_{s_{i+1}}-x_{s_i},\  y_{t_{j+1}}-y_{t_j}\rangle_{\R^d}}{(s_{i+1}-s_i)(t_{j+1}-t_j)}
\end{equation*}
and define the operator $\Lambda^{\le N}_{C_{ij};[s_i,s_{i+1}]\times[t_j,t_{j+1}]}:E_N^2\to E_N^2$ as:
\begin{equation*}
    \Lambda^{\le N}_{C_{ij};[s_i,s_{i+1}]\times[t_j,t_{j+1}]} := \pi_N \circ \Lambda_{C_{ij};[s_i,s_{i+1}]\times[t_j,t_{j+1}]} 
\end{equation*}where $\circ$ denotes composition, $\pi_N$ is the truncation operator (Definition \ref{truncmap}), and $\Lambda_{C_{ij};[s_i,s_{i+1}]\times[t_j,t_{j+1}]}:(\R^2)^\NN\to(\R^2)^\NN$ is the linear operator introduced in Definition \ref{lambda}. Then, for $i\in\{0, \dots, \ell_x-1\}$, $j\in\{0, \dots, \ell_y-1\}$, given initial conditions $\phat^{i,0}=\qhat^{0,j}=(1,0,\dots)\in E_N$ we define recursively 
\begin{equation*}
   (\phat^{i,j+1}, \qhat^{i+1,j}):=\Lambda^{\le N}_{C_{ij};[s_i,s_{i+1}]\times[t_j,t_{j+1}]}(\phat^{i,j}, \qhat^{i,j})
\end{equation*}

where $(\phat^{i,j}, \qhat^{i,j})\equiv (\phat^{i,j}, \qhat^{i,j})^{\le N}\in E_N^2$ by construction. Finally, we define the endpoint of the scheme as:
\begin{equation*}
    \khat^N(s_{\ell_x}, t_{\ell_y}):= \frac{1}{2}\sum_{n=0}^N \left(\phat^{\ell_x-1,\ell_y}_n(s_{\ell_x}-s_{\ell_x-1})^n + \qhat^{\ell_x,\ell_y-1}_n(t_{\ell_y}-t_{\ell_y-1})^n\right)
\end{equation*}
\end{definition}
\begin{remark} Given $\phat^{ij},\ \qhat^{ij}\in E_N$, the map $(\phat^{i,j+1}_n,\qhat^{i+1,j}_n)_{n\in\NN} = \Lambda^{\le N}_{C_{ij};[s_i,s_{i+1}]\times[t_j,t_{j+1}]}(\phat_n^{ij},\ \qhat_n^{ij})_{n\in\NN}$ is explicitly described by the following equations:
\begin{align*}
    \phat^{i,j+1}_n := \sum_{k=0}^{n} \phat^{ij}_k \frac{[C_{ij} (t_{j+1}-t_j)]^{n-k}\ k!}{(n-k)!\ n!} + \sum_{k=1}^{N} \qhat^{ij}_k \frac{ C_{ij}^n(t_{j+1}-t_j)^{n+k}\ k!}{(n+k)!\ n!} \\
    \qhat^{i+1,j}_n := \sum_{k=0}^{n} \qhat^{ij}_k \frac{[C_{ij} (s_{i+1}-s_i)]^{n-k}\ k!}{(n-k)!\ n!} + \sum_{k=1}^{N} \phat^{ij}_k \frac{ C_{ij}^n(s_{i+1}-s_i)^{n+k}\ k!}{(n+k)!\ n!}
\end{align*}
for all $n\le N$, and $\phat^{i,j+1}_n = \qhat^{i+1,j}_n =0$ for $n>N$. 
\end{remark}
\begin{remark}
    Unlike the untruncated version discussed in Section \ref{sec:31}, the symmetry described in Remark \ref{symmetry1} generally does not hold for the polynomial approximation scheme, though we still have: 
\begin{equation*}
 \sum_{n=0}^N \phat^{\ell_x-1,\ell_y}_n(s_{\ell_x}-s_{\ell_x-1})^n \approx \sum_{n=0}^N \qhat^{\ell_x,\ell_y-1}_n(t_{\ell_y}-t_{\ell_y-1})^n
\end{equation*}
The analysis of the associated error estimates is provided later in this section.
\end{remark}
\begin{remark}
    Despite involving the computation of factorials, the scheme can be implemented efficiently by precomputing the factorial terms that are independent of $i$ and $j$. This preprocessing step is independent of the batch-size, the lengths and the dimension of the time series.
\end{remark}

Consider two piecewise-linear $d$-dimensional paths $x$ and $y$  with $\ell_x$ and $\ell_y$ linear segments, respectively. The theoretical time complexity of the polynomial approximation scheme truncated at order $N$ to compute the signature kernel between these two paths is $O(d\ell_x\ell_yN^2)$. Similarly, the finite differences scheme with dyadic refinement, introduced in \cite{salvi2021signature} and outlined here in Definition \ref{dyadic scheme}, has complexity $O(d\ell_x\ell_y2^{2\lambda})$, where $\lambda$ is the dyadic order.

\begin{figure}[!h]
    \centering
    \begin{tikzpicture}[scale=0.8]
        \draw[step=1cm,gray,dotted] (0,0) grid (6,6);
        \draw[step=2cm,black,thin] (0,0) grid (6,6);
        \foreach \i in {0,...,3}
        \foreach \j in {0,...,3}
            \fill (2*\i, 2*\j) circle(.08);
        \foreach \i in {0,...,2}
        \foreach \j in {0,...,2}
        \foreach \s in {0,...,2}
        \foreach \t in {0,...,2}
            \fill (2*\i + \s, 2*\j + \t) circle(.06);
        \draw[thick,->] (0,0) -- (7,0) node[anchor=north west] {$s$};
        \draw[thick,->] (0,0) -- (0,7) node[anchor=south east] {$t$};
        \foreach \i in {0,...,3}{
            \ifnum\i>0{
                \draw[blue, thick]{(0,2*\i) -- (2*\i,0)};
                \draw[orange, thick]{(0,2*\i-1) -- (2*\i-1,0)};
            }
            \fi
            \draw (2*\i,0) node[anchor=north] {$s_{\i}$};
            \draw (0,2*\i) node[anchor=east] {$t_{\i}$};
        }
        \foreach \i in {0,...,3}
        \foreach \j in {0,...,3}
            \fill (10+2*\i, 2*\j) circle(.08);
        \foreach \i in {0,...,3}
        \foreach \j in {0,...,3}
        \draw[step=2cm,black,very thin] (10,0) grid (16,6);
        \draw[thick,->] (10,0) -- (17,0) node[anchor=north west] {$s$};
        \draw[thick,->] (10,0) -- (10,7) node[anchor=south east] {$t$};
        \foreach \i in {0,...,3}{
            \draw (10+2*\i,0) node[anchor=north] {$s_{\i}$};
            \draw (10,2*\i) node[anchor=east] {$t_{\i}$};
            \foreach \j in {0,...,3}{
                \ifnum\j<\i{
                    \ifnum\i=2
                        \draw[blue,very thick]{
                            (10+2*\j,2*\i-2*\j) -- (12+2*\j,2*\i-2*\j)
                            -- (12+2*\j,2*\i-2*\j-2)
                        };
                    \else
                        \draw[orange, very thick]{
                            (10+2*\j,2*\i-2*\j) -- (12+2*\j,2*\i-2*\j)
                            -- (12+2*\j,2*\i-2*\j-2)
                        };
                    \fi
                }
                \fi
            }
        }
    \end{tikzpicture}
    \caption{On the left, the grid illustrates the dyadic approximation scheme for $\lambda = 1$, where the colored lines along the anti-diagonal nodes indicate the nodes that can be computed in parallel when using GPU. On the right, the grid for the polynomial approximation scheme is shown, which enables parallel computation of the coefficients associated with the boundaries (or edges) on the same anti-diagonal. Figure adapted from \cite{jeffrey}.} 
    \label{gridparallel}
\end{figure}
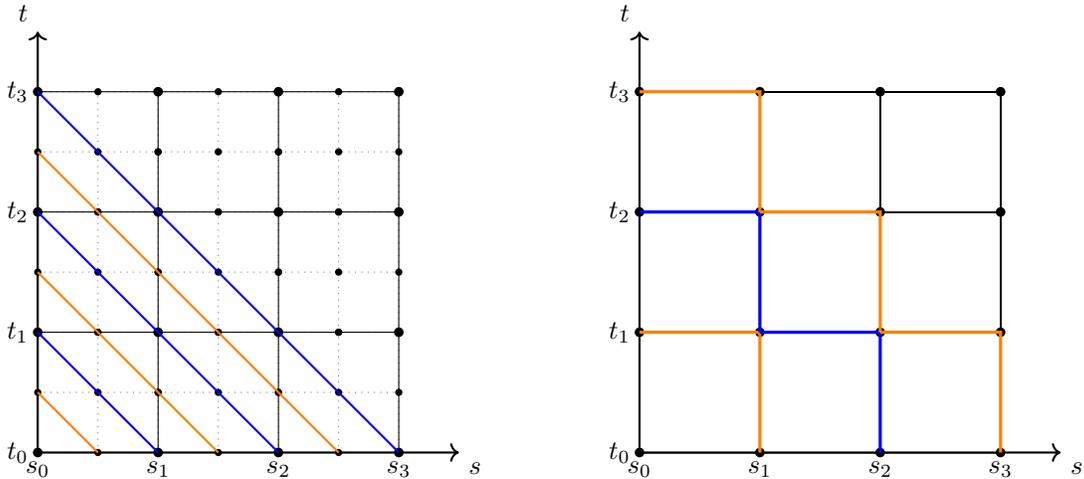

While the finite differences scheme computes the solution at each node, our method computes the coefficients of the polynomials approximating each edge in the grid. Both approaches can leverage parallel computing by updating the anti-diagonals of the solution grid simultaneously. In our method, the polynomial coefficients of the boundaries of the same anti-diagonal cells are updated in parallel (Figure \ref{gridparallel} - right), whereas in the finite differences scheme, the parallelization occurs along the nodes of the same anti-diagonals (Figure \ref{gridparallel} - left). This parallelization breaks the quadratic complexity in the length parameter of both approaches, that becomes linear in $\ell=\max\{\ell_x,\ell_y\}$, making them more practical for large-scale applications. 

\subsubsection*{Error bounds for the polynomial approximation scheme}
In this subsection we define and bound the local and  global truncation error estimate for the polynomial approximation numerical scheme.

The local truncation error quantifies the error introduced in a single approximation step, assuming the initial boundary conditions are exact. More precisely, consider the Goursat PDE on the rectangle $[a,b]\times[c,d]$ as in \eqref{goursatonrectangle} and let $(p,q)\in(\R^2)^\NN$ be the coefficients of the power series describing $k(\cdot,c):[a,b]\to\R$ and $k(a,\cdot):[c,d]\to\R$ respectively. The exact solution at $(b,d)$ is given by
\begin{equation*}
    k(b,d)= \sum_{n=0}^\infty \ptil_n\dab^n = \sum_{n=0}^\infty \qtil_n\dcd^n  
\end{equation*}
with $(\ptil,\qtil)$ from Proposition \ref{approxsolution}. The approximated solution at order $N$, $\khat^N(b,d)$, is here defined as
\begin{equation*}
    \khat^N(b,d):= \frac{1}{2}\sum_{n=0}^N \left(\phat_n\dab^n + \qhat_n\dcd^n\right)
\end{equation*}
where $(\phat, \qhat)= \Lambda^{\le N}_{C,\, [a,b]\times[c,d]}(p,q)^{\le N}$ as in Definition \ref{def:polyapproxscheme}. The local truncation error (LTE) is then defined as
\begin{equation}\label{ltedef}
    LTE:=\left|k(b,d)-\khat^N(b,d)\right|
\end{equation} 

\begin{proposition}[Local Truncation Error]\label{prop:lte}
    Consider the Goursat problem \eqref{goursatonrectangle} on the rectangle $[a,b]\times[c,d]$. Let $(p,q)$ be the coefficients of the power series defining $k(\cdot,c)$ and $k(a,\cdot)$ as in Eq \eqref{initialps} and let $\|\cdot\|_\gamma$ be the norm defined in \eqref{norm}. Assume that $\|(p,q)\|_\gamma\in\R$ for some $0<\gamma\in\N$, and let $\Delta=\max\{b-a, d-c\}$. Then the local truncation error (LTE), defined in \eqref{ltedef}, can be bounded as follows:
    \begin{equation*}
        LTE\le2\|(p,q)\|_\gamma I_0(2\sqrt{\gamma |C|}\Delta)I_0(2\sqrt{(\gamma+1) |C|}\Delta)\frac{(|C|\Delta^2)^{N+1}(\gamma+1)^{N+1}}{[(N+1)!]^2} 
    \end{equation*}
\end{proposition}

\paragraph{Proof.} 

    Fix $2\le N\in\N$ and let $(p,q)\in(\R^2)^\NN$ be given. We define the error coefficients:   
    \begin{equation*}
        (e^1,e^2):=(\ptil,\qtil) - (\phat, \qhat) \in (\R^2)^\NN
    \end{equation*}
    with $(\ptil,\qtil)$ as defined in Proposition \ref{approxsolution}, and $(\phat,\qhat)\equiv(\phat,\qhat)^{\le N}$ representing the approximation of $(\ptil,\qtil)$ obtained by applying the polynomial approximation scheme, i.e. $(\phat,\qhat)=\Lambda^{\le N}_{C,\, [a,b]\times[c,d]}(p,q)^{\le N}$.
    Furthermore, the error terms $(e^1, e^2)$ can be decomposed as:
    \begin{equation*}
        (e^1,e^2)=(e^1,e^2)^{\le N} + (e^1,e^2)^{> N}
    \end{equation*} with
    \begin{align}
        (e^1,e^2)^{\le N} &:=\Lambda^{\le N}_{C,\, [a,b]\times[c,d]}(p,q)^{> N}\equiv (\ptil,\qtil)^{\le N} - (\phat,\qhat) \label{errorcoeff1}\\
        (e^1,e^2)^{> N} &:=\Lambda^{> N}_{C,\, [a,b]\times[c,d]}(p,q)\equiv(\ptil,\qtil)^{>N}\label{errorcoeff2}
    \end{align}
    where $\Lambda^{> N}_{C,\, [a,b]\times[c,d]}=\Lambda_{C,\, [a,b]\times[c,d]}-\Lambda^{\le N}_{C,\, [a,b]\times[c,d]}$.
    
    In summary, $(e^1,e^2)^{\le N}$ accounts for the truncation error arising from truncating the exact initial boundary coefficients $(p,q)$ up to order $N$ - i.e. $(p,q)^{\le N}$, when generating the new coefficients of order $\le N$ - i.e. $(\ptil,\qtil)^{\le N}$. Instead, $(e^1,e^2)^{>N}$ captures the error from truncating the new coefficients $(\ptil, \qtil)$ at level $N$.

Based on the definitions of $(e_1, e_2)$ and local truncation error (Eq. \eqref{ltedef}), the latter can be expressed as the following series:
    \begin{equation*}
        LTE \ =\ \left|\frac{1}{2}\sum_{n=0}^\infty\left( e^1_n\Delta_{a,b}^n +  e^2_n\Delta_{c,d}^n\right) \right|
    \end{equation*} 

To prove the bound on the LTE we will proceed in two steps. First, we calculate the bound for $\|(e^1,e^2)^{\le N}\|_{\gamma+1}$ for a fixed $1\le\gamma\in\N$, then we focus on $\|(e^1,e^2)^{> N}\|_{\gamma+1}$. 

\textbf{\textit{Step 1}}: We aim to bound the norm $\|(e^1,e^2)^{\le N}\|_{\gamma+1}$. For $n\le N$, observe that 
\begin{align*}
    \left|\frac{(n!)^2e^1_n}{(\gamma+1)^n(C\Delta)^n}\right|&\le n!\sum_{k=N+1}^\infty \frac{(k!)^2q_k\Delta_{c,d}^k}{(n+k)!k!(\gamma+1)^n} \\
    &\le n! \sum_{k=N+1}^\infty \frac{|(k!)^2q_k (C\Delta)^{-k}\gamma^{-k}|(\gamma|C|)^k\Delta^{2k}}{n!(k!)^2} \\
    &\le\|(p,q)\|_\gamma\frac{(\gamma|C|\Delta^2)^{N+1}}{[(N+1)!]^2}\sum_{k=0}^\infty \frac{(\gamma|C|)^k\Delta^{2k}}{(k!)^2} \\
    &= \|(p,q)\|_\gamma \frac{(\gamma|C|\Delta^2)^{N+1}}{[(N+1)!]^2} I_0(2\sqrt{\gamma |C|}\Delta)
\end{align*}
A similar argument holds for $e^2_n$. Hence,
\begin{equation}\label{l1te1eq}    
    \|(e^1,e^2)^{\le N}\|_{\gamma+1}\le \|(p,q)\|_\gamma \frac{(\gamma|C|\Delta^2)^{N+1}}{[(N+1)!]^2} I_0(2\sqrt{\gamma |C|}\Delta)
\end{equation}
Thus, we have 
\begin{align}
    \left|\frac{1}{2}\sum_{n=0}^N\left( e^1_n\Delta_{a,b}^n +  e^2_n\Delta_{c,d}^n\right) \right|&\le \sum_{n=0}^N\|(p,q)\|_\gamma I_0(2\sqrt{\gamma |C|}\Delta) \frac{(\gamma|C|\Delta^2)^{N+1}}{[(N+1)!]^2} \frac{(\gamma+1)^n(|C|\Delta)^n}{(n!)^2}\Delta^n \nonumber\\
    & \le\|(p,q)\|_\gamma I_0(2\sqrt{\gamma |C|}\Delta)\frac{(\gamma|C|\Delta^2)^{N+1}}{[(N+1)!]^2}  \sum_{n=0}^\infty\frac{(\gamma+1)^n(|C|\Delta)^n}{(n!)^2}\Delta^n \nonumber \\
    & =\|(p,q)\|_\gamma I_0(2\sqrt{\gamma |C|}\Delta)I_0(2\sqrt{ (\gamma+1)|C|}\Delta)\frac{(\gamma|C|\Delta^2)^{N+1}}{[(N+1)!]^2}  \label{lte1eq}
\end{align}
\textbf{\textit{Step 2}}: Next, consider $(e^1,e^2)^{> N}$, recalling that we have already established the relation $(e^1,e^2)^{> N}=(\ptil,\qtil)^{>N}$ - Eq. \eqref{errorcoeff2}. From Proposition \ref{prop:onestepestimate}, for all $n\ge0$
$$\|(\ptil,\qtil)\|_{\gamma+1} \le \|(p,q)\|_\gamma I_0(2\sqrt{\gamma |C|}\Delta)$$
This implies  
\begin{equation*}
    \max\{\ptil_n,\qtil_n\}\le \|(p,q)\|_\gamma I_0(2\sqrt{\gamma |C|}\Delta)\frac{(\gamma+1)^n(|C|\Delta)^n}{(n!)^2}
\end{equation*}
Therefore, for $n>N$ and $\gamma\ge1$
\begin{align}
    \left|\frac{1}{2}\sum_{n=N+1}^\infty\left( e^1_n\Delta_{a,b}^n +  e^2_n\Delta_{c,d}^n\right) \right|&\le \sum_{n=N+1}^\infty\|(p,q)\|_\gamma I_0(2\sqrt{\gamma |C|}\Delta)\frac{(\gamma+1)^n(|C|\Delta)^n}{(n!)^2}\Delta^n\nonumber\\
    & \le\|(p,q)\|_\gamma I_0(2\sqrt{\gamma |C|}\Delta)  \sum_{n=N+1}^\infty\frac{(\gamma+1)^n(|C|\Delta)^n}{(n!)^2}\Delta^n\nonumber\\
    & \le\|(p,q)\|_\gamma I_0(2\sqrt{\gamma |C|}\Delta)\frac{(\gamma+1)^{N+1}(|C|\Delta^2)^{N+1}}{[(N+1)!]^2}  \sum_{n=0}^\infty\frac{(\gamma+1)^n(|C|\Delta)^n}{(n!)^2}\Delta^n\nonumber\\
    & =\|(p,q)\|_\gamma I_0(2\sqrt{\gamma |C|}\Delta)I_0(2\sqrt{(\gamma+1) |C|}\Delta)\frac{(\gamma+1)^{N+1}(|C|\Delta^2)^{N+1}}{[(N+1)!]^2} \label{lte2eq}
\end{align}
Combining the results \eqref{lte1eq} and \eqref{lte2eq}, by the triangle inequality, we obtain:
$$LTE\le 2\|(p,q)\|_\gamma I_0(2\sqrt{\gamma |C|}\Delta)I_0(2\sqrt{(\gamma+1) |C|}\Delta)\frac{(\gamma+1)^{N+1}(|C|\Delta^2)^{N+1}}{[(N+1)!]^2} $$\qed

To assess the overall performance of the order-$N$ polynomial approximation scheme across the entire grid, we analyze the global truncation error (GTE), which provides a measure of how well this scheme captures the true solution of the signature kernel PDE \eqref{sigpde}. The GTE is defined as:
\begin{equation}\label{gtedef}
    GTE:=\left|k(s_{\ell_x},t_{\ell_y})- \khat^N(s_{\ell_x},t_{\ell_y})\right|
\end{equation} where $k(s_{\ell_x}, t_{\ell_y})$ is the signature kernel for two piecewise linear paths $x$ and $y$ and $\khat^N(s_{\ell_x},t_{\ell_y})$ its polynomial approximation at order $N$ obtained applying the scheme introduced in Definition \ref{def:polyapproxscheme}. 

\begin{theorem}[Global Truncation Error]\label{prop:globalte}

    Let $\khat^N$ be a numerical solution obtained by applying the polynomial approximation scheme truncated at order $N$ (Definition \ref{def:polyapproxscheme}) to the Goursat problem \eqref{sigpde} on $D_I\times D_J$, where the paths $x$ and $y$ are piecewise linear with respect to the grids $\D_I$ and $\D_J$. Recall that $\ell_x=|\D_I|$, $\ell_y=|\D_J|$, and assume there exist a constant $\CC\in \R$ such that 
    \begin{equation*}
        \sup_{\D_I\times\D_J}\left|\left\langle\dot{x}_s, \dot{y}_t\right\rangle\right|\le \CC
    \end{equation*}
    Then 
    \begin{equation*}
        |k(s_{\ell_x},t_{\ell_y})- \khat^N(s_{\ell_x},t_{\ell_y})|\le  f(\ell_x+\ell_y - 1)\frac{(\CC\Delta^2)^{N+1}}{[(N+1)!]^2}\left[\frac{(\ell_x+\ell_y - 1)^{N+2}}{N+2}+(\ell_x+\ell_y - 1)^{N+1}\right]
    \end{equation*}
    where $f$ is defined in Eq. \eqref{fdefinition} and $\Delta=\max_{i,j}\{s_{i+1}-s_i,t_{j+1}-t_j\}$. 
\end{theorem}

\textbf{Proof.} To begin, we present the polynomial approximation scheme using a notation consistent with the one introduced in Definition \ref{lambda} for the exact signature kernel. Let 
\begin{equation*}
        \Lambdahat_{ij} = (\phat^{i,j+1}, \qhat^{i+1,j}) = \Lambda^{\le N}_{C_{ij};[s_i,s_{i+1}]\times[t_j,t_{j+1}]} (\phat^{ij},\qhat^{ij})
\end{equation*} 
where $(\phat^{ij}, \qhat^{ij})\equiv(\phat^{ij}, \qhat^{ij})^{\le N}$ by definition. To simplify notation, $\D_{ij}= [s_i,s_{i+1}]\times[t_j,t_{j+1}]$. Furthermore, let $\CC=\max_{i,j}|C_{ij}|$, as for piecewise linear paths $\sup_{\D_I\times\D_J}\left|\left\langle\dot{x}_s, \dot{y}_t\right\rangle\right|=\max_{i,j}|C_{ij}|$. The error coefficients are introduced as the differences between the exact and approximated coefficients, and are defined as follows:
\begin{equation*}
    \Xi_{ij} = (\mu^{i,j+1},\nu^{i+1,j}):=\Lambda_{ij} - \Lambdahat_{ij}= (p^{i,j+1}-\phat^{i,j+1},q^{i+1,j}-\qhat^{i+1,j})\in (\R^2)^\NN
\end{equation*}
Notice that by construction, since $(p^{ij}, q^{ij})\in E_N\times E_N$ - where $E_N=\{e\in\R^\NN: e=e^{\le N}\}$ - it follows that $(\mu^{ij}, \nu^{ij})^{>N}=(p^{ij}, q^{ij})^{>N}$. Then by the definition of the global truncation error (Eq. \eqref{gtedef}), and the polynomial approximation scheme (Definition \ref{def:polyapproxscheme}):
\begin{equation*}
    \left|k(s_{\ell_x},t_{\ell_y})- \khat^N(s_{\ell_x},t_{\ell_y})\right|=\left|\frac{1}{2}\sum_{n=0}^\infty\left( \mu^{l_x-1,l_y}_n(s_{l_x}-s_{l_x-1})^n +  \nu^{l_x,l_y-1}_n(t_{l_y}-t_{l_y-1})^n \right)\right|
\end{equation*}
This can be decomposed into two terms using the triangle inequality, leading to the bound $GTE\le GTE_1+GTE_2$, where:
\begin{equation*}
    GTE_1:=\left|\frac{1}{2}\sum_{n=0}^N\left( \mu^{l_x-1,l_y}_n(s_{l_x}-s_{l_x-1})^n +  \nu^{l_x,l_y-1}_n(t_{l_y}-t_{l_y-1})^n \right)\right|
\end{equation*}
and
\begin{equation*}
    GTE_2:=\left|\frac{1}{2}\sum_{n=N+1}^\infty\left( \mu^{l_x-1,l_y}_n(s_{l_x}-s_{l_x-1})^n +  \nu^{l_x,l_y-1}_n(t_{l_y}-t_{l_y-1})^n \right)\right|
\end{equation*}
Now we are interested in analysing the coefficients $(\mu^{ij}, \nu^{ij})\in(\R^2)^\NN$ for all $i,j$ for which are defined. Notice that for all $i=0,1,\dots,\ell_x-1, j=0,1,\dots,\ell_y-1$, it holds that
\begin{equation}\label{xxi}
\Xi_{ij} = \Lambda^{\le N}_{C_{ij};\D_{ij}}(\mu^{ij},\nu^{ij})^{\le N} + \Lambda^{\le N}_{C_{ij};\D_{ij}}(\mu^{ij},\nu^{ij})^{> N} + \Lambda^{> N}_{C_{ij};\D_{ij}}(p^{ij},q^{ij}) 
\end{equation}
where $\Lambda^{> N}_{C_{ij};\D_{ij}}:= \Lambda_{C_{ij};\D_{ij}} - \Lambda^{\le N}_{C_{ij};\D_{ij}}$.

Equation \eqref{xxi} means that at each step in the approximation scheme three types of error contribute to $\Xi_{ij}=(\mu^{i,j+1}, \nu^{i+1,j})$:

\begin{enumerate}
    \item Error propagation from previous steps: this affects coefficients up to order N and is represented as $$\Lambda^{\le N}_{C_{ij};\D_{ij}}
    (\mu^{ij},\nu^{ij})^{\le N}$$ 
    \item Error from truncating initial coefficients: this error arises within the current rectangle due to truncation of the exact initial boundary coefficients $(p^{ij}, q^{ij})$ at level $N$. It affects the terms of $\Xi_{ij}$ up to order $N$ and is represented by: $$\Lambda^{\le N}_{C_{ij};\D_{ij}}(\mu^{ij},\nu^{ij})^{> N}\equiv \Lambda^{\le N}_{C_{ij};\D_{ij}}(p^{ij},q^{ij})^{> N}$$ 
    as by construction $(p^{ij},q^{ij})^{> N}\equiv (\mu^{ij},\nu^{ij})^{> N}$.
    \item Error from truncating the exact solution beyond order $N$: this captures errors from omitting higher-order terms and is represented as $$\Lambda^{> N}_{C_{ij};\D_{ij}}(p^{ij},q^{ij})=(p^{i,j+1},q^{i+1,j})^{>N}$$ 
\end{enumerate}

Notice that only the first two terms in the RHS of \eqref{xxi} contribute to $\Xi_{ij}^{\le N}$. To better understand the behaviour of these first two terms, first let 
\begin{equation*}
    \xi:= \frac{(\CC\Delta^2)^{N+1}}{[(N+1)!]^2} 
\end{equation*}
where recall that $K=\max_{i,j}|C_{i,j}|$. Then,
\begin{align*}
    \|\Lambda^{\le N}_{C_{ij};\D_{ij}}(\mu^{ij},\nu^{ij})^{\le N}\|_{i+j+1} &\le \|(\mu^{ij},\nu^{ij})^{\le N}\|_{i+j}I_0(2\sqrt{(i+j)\CC}\Delta) \\
    \|\Lambda^{\le N}_{C_{ij};\D_{ij}}(p^{ij},q^{ij})^{> N}\|_{i+j+1} &\le \|(p^{ij},q^{ij})\|_{i+j}\ \xi I_0(2\sqrt{(i+j)\CC}\Delta)\ (i+j)^{N+1}
\end{align*}
where the first line follows from Proposition \ref{prop:onestepestimate}, which also applies to the coefficients of the approximation scheme since $\|\Lambda^{\le N}_{C_{ij};\D_{ij}}\|_{op}\le \|\Lambda_{C_{ij};\D_{ij}}\|_{op}$, and the second line follows from the proof of the local truncation error - Eqs. \eqref{errorcoeff1} and \eqref{l1te1eq}. 
By the triangle inequality:
\begin{align*}
    \|\Xi_{ij}^{\le N}\|_{i+j+1} &\le    \|(\mu^{ij},\nu^{ij})^{\le N}\|_{i+j}I_0(2\sqrt{(i+j)\CC}\Delta) + \|(p^{ij},q^{ij})\|_{i+j}\ \xi I_0(2\sqrt{(i+j)\CC}\Delta)\ (i+j)^{N+1} \\
    &\le \|(\mu^{ij},\nu^{ij})^{\le N}\|_{i+j}I_0(2\sqrt{(i+j)\CC}\Delta) + \xi f(i+j)\ (i+j)^{N+1}
\end{align*}
with $f$ defined in \eqref{fdefinition} and where we have used the result $\|(p^{ij},q^{ij})\|_{i+j}\le f(i+j-1)$ from Proposition \ref{multistepestimate}. Moreover, a straightforward computation shows that $\|\Xi_{00}^{<N}\|_1=0$.

Using a similar induction to the one in the proof of Proposition \ref{multistepestimate}, we get
\begin{equation*}
    \|\Xi_{ij}^{\le N}\|_{i+j+1}\le \xi f(i+j) \sum_{m=0}^{i+j} m^{N+1}\le \xi f(i+j)\ \frac{(i+j+1)^{N+2}-1}{N+2}
\end{equation*}
This leads to the following bound:
\begin{align*}
    GTE_1&\le\left|\frac{1}{2}\sum_{n=0}^N\left( \mu^{l_x-1,l_y}_n\Delta^n +  \nu^{l_x,l_y-1}_n\Delta^n \right)\right| \\
    &\le f(l_x+l_y-2)\xi\ \frac{(l_x+l_y-1)^{N+2}-1}{N+2}{}\sum_{n=0}^N\frac{(l_x+l_y-1)^n(\CC\Delta)^n\Delta^n}{(n!)^2} \\
    &\le f(l_x+l_y-1)\frac{(l_x+l_y-1)^{N+2}(\CC\Delta^2)^{N+1}}{(N+2)[(N+1)!]^2}
\end{align*}

On the other hand,  we have established that $\Xi_{ij}^{>N}=\Lambda_{ij}^{>N}=(p^{ij}, q^{ij})^{>N}$. From Proposition \ref{multistepestimate}, we have 
\begin{align*}
\|\Xi_{ij}^{>N}\|_{i+j+1}&=\|\Lambda_{ij}^{>N}\|_{i+j+1}\le \|\Lambda_{ij}\|_{i+j+1} \le f(i+j)
\end{align*}
Thus, for $i=l_x-1$ and $j=l_y-1$ 
\begin{align*}
    GTE_2&\le\left|\frac{1}{2}\sum_{n=N+1}^\infty\left( \mu^{l_x-1,l_y}_n\Delta^n +  \nu^{l_x,l_y-1}_n\Delta^n \right)\right| \\
    &\le f(l_x+l_y-2)\sum_{n=N+1}^\infty\frac{(l_x+l_y-1)^n(\CC\Delta)^n\Delta^n}{(n!)^2} \\
    &\le f(l_x+l_y-1)\frac{(l_x+l_y-1)^{N+1}(\CC\Delta^2)^{N+1}}{[(N+1)!]^2}
\end{align*}
Then, as $GTE\le GTE_1 + GTE_2$, we obtain:
\begin{align*}
    GTE \le f(l_x+l_y - 1)\frac{(\CC\Delta^2)^{N+1}}{[(N+1)!]^2}\left[\frac{(l_x+l_y - 1)^{N+2}}{N+2}+(l_x+l_y - 1)^{N+1}\right]
\end{align*}\qed

Theorem \ref{prop:globalte} provides the convergence result for our numerical scheme. For a fixed constant $K\in\R$, taking the limit as $N\to\infty$  yields convergence to the true solution. In Section \ref{sec:exp}, we empirically demonstrate that increasing the hyperparameter $N$ in the polynomial approximation scheme leads to faster convergence compared to increasing the dyadic order (Definition \ref{dyadic scheme}) in the finite differences scheme. Moreover, the finite differences scheme and the polynomial approximation scheme are not mutually exclusive; they can be combined to produce a hybrid method by introducing a refined grid within the polynomial approximation scheme.

\subsection{Polynomial interpolation numerical scheme}\label{sec:polyinterp}

In Proposition \ref{approxsolution}, we derived an analytical representation of the solution to the Goursat problem on a rectangular domain - Eq. \eqref{goursatonrectangle},  which involved expressing the initial boundary conditions as power series. By truncating these series, we constructed the polynomial approximation scheme - Definition \ref{def:polyapproxscheme}, which outperforms finite differences methods.

Alternatively, when the initial boundary conditions are polynomials, the exact solution on the opposite boundaries can be expressed in terms of special functions, specifically the confluent hypergeometric limit function \cite{andrews1999special}.

\begin{proposition}\label{propinterp}
Consider the Goursat problem on a rectangle $[a,b]\times[c,d]$ as in \eqref{goursatonrectangle}. Further, assume that the initial boundary conditions $k(\cdot,c):[a,b]\to\R$ and $k(a,\cdot):[c,d]\to\R$ are polynomials at $s=a$ and $t=c$ respectively, such that
\begin{equation*}
    k(s,c)=\sum_{n=0}^N g_n (s-a)^n  \qquad k(a,t) =\sum_{n=0}^N h_n (t-c)^n
\end{equation*}with the constraint $g_0=h_0$. Then the solution $k(s,t)$ can be written as
\begin{equation}\label{interpformula}
    k(s,t) = g_0I_0(2\sqrt{C(s-a)(t-c)}) + \sum_{n=1}^{N} \left(g_n (s-a)^n + h_n(t-c)^n \right)\hypf(n+1;\, C(s-a)(t-c)) 
\end{equation}
with $s\in[a,b]$, $t\in[c,d]$. Here, $\hypf$ is the confluent hypergeometric limit function \cite{andrews1999special} defined by:
\begin{equation*}
    \hypf(k;z)=\sum_{n=0}^\infty\frac{z^n}{(k)_nn!}
\end{equation*} where $(k)_n$ is the Pochhammer symbol, which represents the rising factorial.
\end{proposition}

\begin{remark}
    From the definition of the confluent hypergeometric limit function and the modified Bessel function (Remark \ref{mbf}), we observe that for $z \geq 0$, the identity $\hypf(1, z) = I_0(2\sqrt{z})$ holds.
\end{remark}
\begin{remark}
    See Chapter 3 of the Ph.D Thesis of the third-name author \cite{jeffrey}.
\end{remark}

The proof of Proposition \ref{propinterp} can be found in Appendix \ref{app:A1}. 

Moreover, note that for bounded coefficients $(g_n,h_n)_{n=0,\dots,N}$, the function described by \eqref{interpformula} is infinitely differentiable. To streamline the forthcoming analysis, we introduce specific notation to explicitly represent the map $(k(\cdot,c),\,k(a,\cdot))\mapsto(k(\cdot,d),\,k(b,\cdot))$, both in the case where the inputs are polynomials (as outlined in Proposition \ref{propinterp}) and when they are differentiable functions, as described in Theorem \ref{besselsolution}.

\begin{definition}\label{operatorphi}
    Consider the Goursat problem on a general rectangular domain $[a,b]\times[c,d]$ with associated constant $C\in\R/\{0\}$ as in Eq. \eqref{goursatonrectangle}. Let \(\sigma := k(\cdot,c)\in C^1([a,b],\R)\) and \(\tau := k(a,\cdot)\in C^1([c,d],\R)\) denote the boundary conditions at the bottom and left edges of the rectangle. Let \(\sigmanew := k(\cdot,d)\in C^1([a,b],\R)\) and \(\taunew := k(b,\cdot)\in C^1([c,d],\R)\) represent the corresponding solutions at the top and right edges, respectively.
    
    We define the operator $\Phi_{C,[a,b]\times[c,d]}: C^1([a,b],\R)\times C^1([c,d],\R)\to C^1([a,b],\R)\times C^1([c,d],\R)$ such that:
    \begin{equation*}
        (\sigmanew,\taunew)=\Phi_{C,[a,b]\times[c,d]}(\sigma, \tau)
    \end{equation*}
    Similarly we define the restriction of $\Phi_{C;[a,b]\times[c,d]}$ to the space of polynomials of degree at most $N$ with real coefficients, denoted as $P_N(\R)$, as \begin{equation*}
        \Phi^N_{C;[a,b]\times[c,d]}:P_N(\R)\times P_N(\R)\to C^\infty([a,b],\R)\times C^\infty([c,d],\R)
    \end{equation*} which is constructed by applying Proposition \ref{propinterp} to the polynomial initial data.

This construction and notation can be naturally extended to any rectangle $[s_i,s_{i+1}]\times[t_j,t_{j+1}]$ in a given grid, with an associated constant $C_{ij}$, to solve the Goursat PDE on a partitioned domain.
\end{definition}

\begin{remark}
    When the initial boundary conditions $\sigma$ and $\tau$ can be expressed in the form of power series, the map $\Phi_{C; [a,b]\times[c,d]}$ is similar to the map $\Lambda_{C; [a,b]\times[c,d]}(\R^2)^\NN\to(\R^2)^\NN$. The key distinction lies in their domains and outputs: while $\Phi$ and $\Phi^N$ operates on the input functions directly, $\Lambda$ acts on the coefficients of their power series and outputs the coefficients of the power series describing the solution.
\end{remark}

This setup allows us to define an alternative numerical scheme based on polynomial interpolation. 
More precisely, on each rectangle, the polynomial interpolation scheme operates by first constructing the interpolating polynomial of the initial boundaries, and then by using Proposition \ref{propinterp} to evaluate the solution on the opposite boundaries.

To construct the interpolating polynomial, we use the $N$-Chebyshev nodes of the second kind, also called \textit{Chebyshev extrema} \cite{trefethen2019approximation}, which, when scaled to a compact interval $[u, v]$, are defined as
\begin{equation}\label{chebextr}
    \omega_{m;[u,v]} = \frac{v-u}{2}\cos(\pi m/N)+\frac{u+v}{2},\qquad m=0,1,\dots,N-1
\end{equation}

This set of points minimizes the problem of Runge’s phenomenon \cite{runge1901empirische}, which is the large oscillation that can occur at the edges of an interval when using polynomial interpolation. Chebyshev nodes thus provide a near-optimal distribution, reducing interpolation error and enhancing the stability of the polynomial approximation. For convenience, we define a general interpolation operator based on the Chebyshev extrema.

\begin{definition}\label{interpooperator}
     Define $\Omega^N_{[a,b]}:C([a,b],\R)\to P_N(\R)$ as a general interpolation operator that, given a continuous function $f\in C([a,b],\R)$,  evaluates $f$ at the $N+1$ Chebyshev extrema rescaled to the interval $[a,b]$, i.e. $\{\omega_{m,[a,b]}\}_{m=0,\dots,N}$ as defined in \eqref{chebextr}, and returns the unique polynomial of degree at most $N$ that interpolates these points. 
\end{definition}

Building on these definitions and Proposition \ref{propinterp} we define the polynomial interpolation scheme for the signature kernel PDE \eqref{sigpde}. 

\begin{definition}[Polynomial interpolation numerical scheme]\label{def:polyinterpscheme}
Let $\D_I =\{s_0<s_1<\dots<s_{\ell_x}\}$ and $\D_J =\{t_0<t_1<\dots<t_{\ell_y}\}$ be partitions of the compact intervals $I$ and $J$ and let $x:I\to\R^d$ and $y:J\to\R^d$ be two piecewise linear functions on $\D_I$ and $\D_J$ respectively. 
Here, let $\khat^N$ be the solution of the order-$N$ polynomial interpolation scheme to the Goursat PDE \eqref{sigpde} on the grid $\D_I\times\D_J$. 

On each rectangular subdomain \(\D_{ij}=[s_i,s_{i+1}]\times[t_j,t_{j+1}]\), for a fixed $N\ge 2$, we define:
\begin{equation*}
    \sigmahat_{ij} := \khat^N(\cdot, t_j) : [s_i,s_{i+1}]\to\R\qquad\textrm{and}\qquad \tauhat_{ij} := \khat^N(s_i, \cdot) : [t_j,t_{j+1}]\to\R
\end{equation*} Furthermore, let
\begin{equation*}
    C_{ij} = \frac{\langle x_{s_{i+1}}-x_{s_i},\  y_{t_{j+1}}-y_{t_j}\rangle_{\R^d}}{(s_{i+1}-s_i)(t_{j+1}-t_j)}
\end{equation*}
Starting from the initial conditions \(\sigmahat_{i,0}=1\) and \(\tauhat_{0,j}=1\), we proceed recursively as follows: for $i\in\{0, \dots, \ell_x-1\}$, $j\in\{0, \dots, \ell_y-1\}$
\begin{align*}
    (\sigmahatch_{ij},\tauhatch_{ij})&=(\Omega^N_{[s_i,s_{i+1}]}(\sigmahat_{i,j}),\ \Omega^N_{[t_j,t_{j+1}]}(\tauhat_{i,j})) \\
    (\sigmahat_{i,j+1},\tauhat_{i+1,j}) &= \Phi^N_{C_{ij};[s_i,s_{i+1}]\times[t_j,t_{j+1}]}(\sigmahatch_{ij}, \tauhatch_{ij}) 
\end{align*}
where $\Phi^N$ is introduced in Definition \ref{operatorphi}, and $\Omega^N$ is the interpolation operator introduced in Definition \ref{interpooperator}. Consequently,  $\sigmahatch_{ij},\tauhatch_{ij}\in P_N(\R)$ are the polynomial interpolants of $\sigmahat_{ij}$ and $\tauhat_{ij}$, respectively. 
The endpoint of our scheme is then defined as:
\begin{equation*}
    k(s_{\ell_x},t_{\ell_y}):=\frac{1}{2}\left(\sigmahat_{\ell_x-1,\ell_y}(s_{\ell_x})+\tauhat_{\ell_x,\ell_y-1}(t_{\ell_y})\right)
\end{equation*}
Computationally, the polynomial interpolants $\sigmahatch_{ij},\tauhatch_{ij}\in P_N(\R)$ are stored by their coefficients, while the evaluation of $\sigmahat_{ij}\in C^\infty([s_i,s_{i+1}],\R)$ and $\tauhat_{ij}\in C^\infty([t_j,t_{j+1}],\R)$ is required only at the scaled Chebyshev extrema for interpolation on each respective interval.
\end{definition}

\section{Experiments}\label{sec:exp}

In this section, we present the results of our numerical schemes and compare them to the finite differences scheme proposed by \cite{salvi2021signature} and implemented in the \texttt{sigkerax} Python library \cite{sigkerax} - the JAX version of the original Pytorch-based \texttt{sigkernel} library \cite{sigkernel}.

\begin{remark}
    The finite difference scheme implemented in the \texttt{sigkerax} package \cite{sigkerax} introduces a modification from the original scheme (Definition \ref{dyadic scheme}) by replacing the dyadic refinement hyperparameter with a more general refinement factor. Specifically, while the dyadic hyperparameter $\lambda$ increases the grid size by a factor of $2^{2\lambda}$, in the modified algorithm, the \textit{refinement factor} $\gamma$ interpolates each time series by adding a number of points to the grid proportional to $\gamma^2$.
\end{remark}

Additionally, the finite differences implementation in the \texttt{sigkerax} package employs the second-order explicit finite differences scheme detailed in Eq. \eqref{secondexplicitscheme}.

\subsubsection*{Python library}
The code for this paper generated a Python library named Polysigkerax, which provides all the necessary tools to utilize the signature kernel. The code is freely available for download at the link \url{https://github.com/FrancescoPiatti/polysigkernel}. The implementation is written in JAX, enabling multi-threading and multi-processing for large batches. Similarly to the original library, additional features to perform statistical tests and compute mean maximum discrepancy have also been incorporated.

All the experiments have been run using a MacBook Air 2022 with M2 chip (for CPU) and a NVIDIA GTX 3090 (for GPU).

\subsubsection*{Error analysis and computational performance analysis}
We first compare the numerical error of our schemes against the benchmark, which is the truncated signature kernel at level 21. The factorial decay of the signature ensures machine-level precision at this truncation, offering a highly accurate benchmark. However, due to the curse of dimensionality of the signature transform, this benchmark can only be computed for low-dimensional time series, hence our analysis focuses on 2-$d$ paths. Since the dimensionality of the paths enters the algorithms only through the (weighted) dot product of the derivatives, we do not expect higher dimensions to significantly impact the performance of our scheme.

\begin{figure}[!h]
    \centering
    \includegraphics[width=0.9\textwidth]{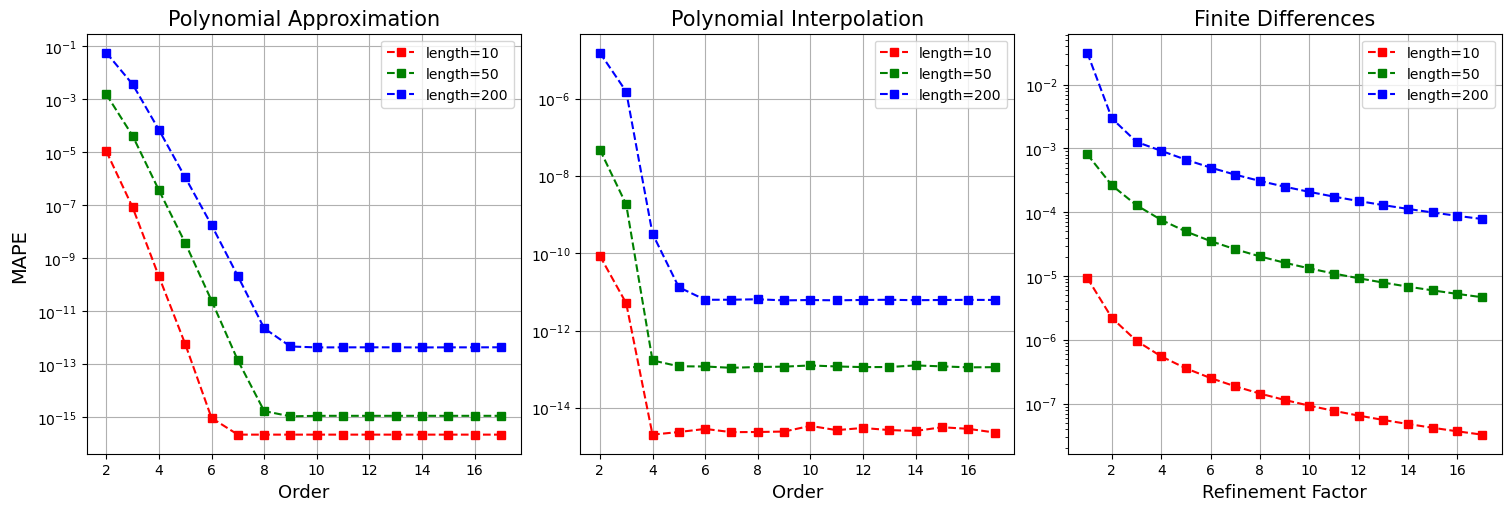}
    \caption{Comparison of MAPE against polynomial order for the polynomial approximation and interpolation schemes is shown on the \textbf{left} and \textbf{center}, respectively. On the \textbf{right}, MAPE is plotted against the refinement factor for the finite differences scheme. This analysis is performed on a $8\times8$ batch of 2-$d$ standard Brownian motions with different lengths (color-coded).}
    \label{fig:errorvslevel}
    \vspace{-0.1in}
\end{figure}

Figure \ref{fig:errorvslevel} illustrates the decrease in the Mean Absolute Percentage Error (MAPE) for both the polynomial approximation and interpolation schemes as the polynomial order increases, compared to the finite differences scheme against growing refinement factors. For shorter time series, our algorithm easily achieves machine precision (approximately $5 \times 10^{-16}$), making it $\approx10^8$ times more accurate than the finite differences scheme, even at high refinement levels. For longer time series, our algorithm maintains robust performance, with the slightly larger errors likely due to the numerical inaccuracies of the benchmark given the extensive computations involved. Tables \ref{tab1}, \ref{tab2} and \ref{tab3} in Appendix \ref{addexp} provide the numerical values for this experiment. Additionally, Appendix \ref{addexp} includes results from a complementary experiment conducted on stochastic paths other than Brownian motions.

As mentioned in Section \ref{sec:polyapprox}, the computational performance of the polynomial approximation scheme is comparable to that of the finite differences scheme. In contrast, the complexity of the polynomial interpolation scheme depends heavily on the choice of the interpolation function $\Omega^N$ (Definition \ref{interpooperator}). In our experiments, we use \texttt{jax.numpy.polyfit} \cite{jax2018github}, which computes the monomial basis coefficients via a least-squares problem. This choice, combined with our implementation of the hypergeometric function \(\hypf\), introduces additional computational complexity as the interpolation order increases, rendering it slower compared to the other two schemes. While more efficient implementations of $\Omega$ and \(\hypf\) could mitigate this issue, the computational time analysis presented below is focused exclusively on the polynomial approximation scheme.

\begin{figure}[!h]
    \centering
    \includegraphics[width=0.9\textwidth]{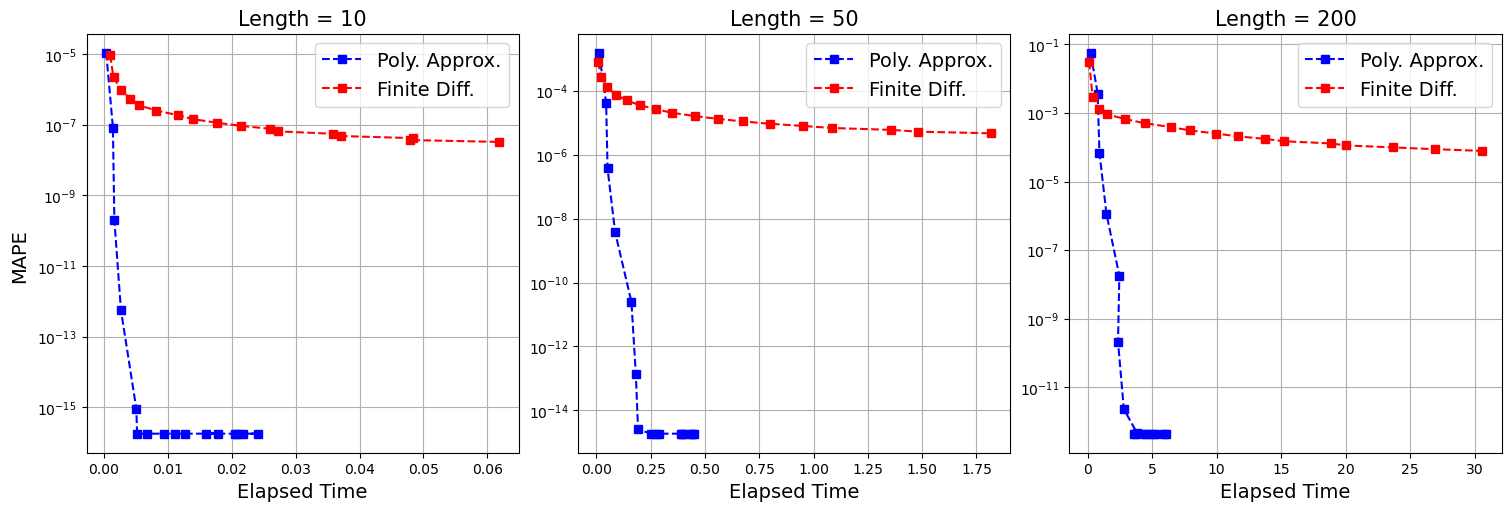}
    \caption{Mean Absolute Percentage Error (MAPE) against time (s) for increasing polynomial orders and refinement factors in the polynomial approximation scheme and the finite differences scheme with grid refinement, respectively. The three plots display the results for realizations of a $4 \times 4$ batch of 2-dimensional Brownian motions across different lengths $\ell$: $\ell=10$ (\textbf{left}), $\ell=50$ (\textbf{middle}), and $\ell=200$ (\textbf{right}).}
    \label{fig:errorvstime}
\end{figure}

The key result of our experiments is that the polynomial approximation scheme achieves significantly higher accuracy while being more computationally efficient than the finite differences scheme. Figure \ref{fig:errorvstime} shows the Mean Absolute Percentage Error (MAPE) against computational time for both schemes: the polynomial approximation method with increasing polynomial orders and the finite differences method with increasing refinement factors. 

To further investigate the complexities of the two algorithms, we distinguish between theoretical and empirical complexities. Theoretically, the polynomial approximation scheme has a complexity of  $O(d\ell_x\ell_yN^2)$, where $\ell_x$ and $\ell_y$ are the lengths of the two time series, while the finite differences scheme has complexity $O(d\ell_x\ell_y\gamma^2)$. Both schemes exhibit quadratic complexity with respect to their hyperparameters, $N$ and $\gamma$. However, these theoretical estimates are asymptotic and do not always reflect practical performance, especially when the hyperparameters are relatively small. Empirically, our experiments reveal that for smaller datasets, the polynomial approximation scheme exhibits near-linear complexity in $N$, giving our algorithm a significant performance edge over the finite differences scheme on the CPU. This is illustrated in Figure \ref{fig:timeanal2} and further discussed in Appendix \ref{addexp}, where a more detailed analysis of the computational time for our scheme is also presented.

Furthermore, in the context of GPU acceleration, enabling the parallelization described in Section \ref{sec:polyapprox} and Figure \ref{gridparallel} breaks the quadratic complexity in relation to the time series lengths, reducing it to linear complexity (see Appendix \ref{addexp}).

\section*{Acknowledgements}
TC has been supported by the EPSRC Programme Grant EP/S026347/1 and acknowledges the support of
the Erik Ellentuck Fellowship at the Institute for Advanced Study, where part of this work was carried out. FP has been supported by the EPSRC Centre for Doctoral Training in Mathematics of Random Systems: Analysis, Modelling and Simulation (EP/S023925/1). JP has been supported by the EPSRC CDT in Financial Computing and Analytics (EP/L015129/1). For the purpose of open access, the authors have applied a Creative Commons Attribution (CC BY) licence to any Author Accepted Manuscript version arising.

\printbibliography[heading=bibintoc,title={References}]

\clearpage
\appendix
\setstretch{1.25}
\section{Proofs}\label{app:A1}
In this section, we provide the proofs for the key results presented in this paper, specifically Proposition \ref{approxsolution} and \ref{propinterp}, where we derive exact solutions to the Goursat problem on a rectangle for piecewise linear time series.

We begin by stating and proving Proposition \ref{prop:a1}, which serves as a foundational result. 
\begin{proposition}\label{prop:a1}
    Consider the Goursat problem \eqref{goursatonrectangle} on the rectangular domain $\mathcal{D}_R=\left\{(s, t) \mid a \leq s \leq b, c \leq t \leq d\right\}$. Assume that $\sigma$ and $\tau$ are (real) analytic functions at $s=a$, $t=c$ respectively and are determined by their power series
    \begin{equation}\label{appeq1}
        \sigma(s) = k(s, c) = \sum_{n\in\NN} p_n (s-a)^n \qquad \tau(t) = k(a, t) = \sum_{n\in\NN} q_n (t-c)^n
    \end{equation}
    with the constrain $\sigma(a)=\tau(c) \implies p_0=q_0$, as $\sigma(a)=\tau(c)=k(a,c)$. Furthermore, assume that these power series converge uniformly on $[a,b]$ and $[c,d]$ respectively. Then
    \begin{align}\label{appeq2}
        k(s,t)&= \sum_{n=1}^{\infty} \sum_{m=0}^{\infty} p_n \frac{[C(t-c)]^{m} (s-a)^{n+m}}{m!} \frac{n!}{(n+m)!} + \sum_{n=0}^{\infty} \sum_{m=0}^{\infty} q_n \frac{[C (s-a)]^{m} (t-c)^{n+m}}{m!} \frac{n!}{(n+m)!}
    \end{align} 
    or, equivalently,
    \begin{align}\label{appeq3}
        k(s,t)&= \sum_{n=0}^{\infty} \sum_{m=0}^{\infty} p_n \frac{[C(t-c)]^{m} (s-a)^{n+m}}{m!} \frac{n!}{(n+m)!} + \sum_{n=1}^{\infty} \sum_{m=0}^{\infty} q_n \frac{[C (s-a)]^{m} (t-c)^{n+m}}{m!} \frac{n!}{(n+m)!}
    \end{align} 
\end{proposition}

\paragraph{Proof.} The proof is organised in two steps. In the first step we demonstrate that 
\begin{equation}\label{step1eq1}
        k(s,t)=-\sigma(a) I_0\left(2 \sqrt{C(s-a) (t-c)} \right) 
        + \sigma(s) + \tau(t) + Q_1(s,t) + Q_2(s,t)
    \end{equation} where 
    \begin{align*}
        Q_1(s,t)&=2\sqrt{C(s-a)(t-c)}\intopi \sigma((s-a)\cos^2\theta+a)I_1\left(2 \sqrt{C(s-a) (t-c)}\sin\theta\right)\cos\theta d\theta \\
        Q_2(s,t)&=2\sqrt{C(s-a)(t-c)}\intopi \tau((t-c)\cos^2\theta+c)I_1\left(2 \sqrt{C(s-a) (t-c)}\sin\theta\right)\cos\theta d\theta
    \end{align*}
In the second step we substitute the power series \eqref{appeq1} for the initial boundaries into Eq. \eqref{step1eq1}, recovering \eqref{appeq2} and \eqref{appeq3}.

\textbf{\textit{Step 1}}: From Theorem \ref{besselsolution}, for $s\in[a,b]$, $t\in[c,d]$, we have 
\begin{align*}
k(s, t)&=\sigma(a) I_0\left(2 \sqrt{C(s-a)(t-c)} \right) +\int_a^s \sigma^{\prime}(r) I_0\left(2 \sqrt{C(s-r)(t-c)}\right) d r +\int_c^t \tau^{\prime} (r) I_0 \left(2 \sqrt{C(s-a)(t-r)}\right) d r 
\end{align*}

where $\sigma(\cdot)=k(\cdot,c)$ and $\tau(\cdot)=k(a,\cdot)$. Integrating by parts, we obtain
\begin{align*}
    k(s,t)&=\sigma(a) I_0\left(2 \sqrt{C(s-a) (t-c)} \right) + \sigma(s)I_0(0) + \tau(t)I_0(0)
    \\&- \sigma(a) I_0\left(2 \sqrt{C(s-a) (t-c)} \right) - \tau(c) I_0\left(2 \sqrt{C(s-a) (t-c)} \right) 
    \\ &+\sqrt{C(t-c)}\int_a^s \frac{\sigma(r)}{\sqrt{s-r}} I_1\left(2 \sqrt{C(s-r) (t-c)}\right)d r +\sqrt{C(s-a)}\int_c^t \frac{\tau(r)}{\sqrt{t-r}} I_1 \left(2 \sqrt{C(s-a)(t-r)}\right)d r
\end{align*}

Given that $\tau(c) = \sigma(a)$ and recalling that $I_0(0)=1$, this simplifies to
\begin{align}
    k(s,t)&=-\sigma(a) I_0\left(2 \sqrt{C(s-a) (t-c)} \right) 
    + \sigma(s) + \tau(t)  \label{integralsol}  \\ &+\sqrt{C(t-c)}\int_a^s \frac{\sigma(r)}{\sqrt{s-r}} I_1\left(2 \sqrt{C(s-r) (t-c)}\right)d r + \sqrt{C(s-a)}\int_c^t \frac{\tau(r)}{\sqrt{t-r}} I_1 \left(2 \sqrt{C(s-a)(t-r)}\right)d r \nonumber
\end{align}

Now define:
\begin{equation*}
Q_1(s,t):=\sqrt{C(t-c)}\int_a^s \frac{\sigma(r)}{\sqrt{s-r}} I_1\left(2 \sqrt{C(s-r) (t-c)}\right)d r\end{equation*}
\begin{equation*}
Q_2(s,t) := \sqrt{C(s-a)}\int_c^t \frac{\tau (r)}{\sqrt{t-r}} I_1 \left(2 \sqrt{C(s-a)(t-r)}\right)d r\end{equation*}
Using the change of variables:
\begin{align*}
\sqrt{s-a} \sin \theta=\sqrt{s-r} \qquad\textrm{and}\qquad \sqrt{s-a} \cos \theta=\sqrt{r-a}
\end{align*}
we obtain
\begin{equation*}
Q_1(s,t)=2\sqrt{C(s-a)(t-c)}\intopi \sigma((s-a)\cos^2\theta+a)I_1\left(2 \sqrt{C(s-a) (t-c)}\sin\theta\right)\cos\theta d\theta
\end{equation*}
A similar change of variables for $Q_2(s)$ yields:
\begin{equation*}Q_2(s,t)=2\sqrt{C(s-a)(t-c)}\intopi \tau((t-c)\cos^2\theta+c)I_1\left(2 \sqrt{C(s-a) (t-c)}\sin\theta\right)\cos\theta d\theta\end{equation*}

\textbf{\textit{Step 2:}} Substituting the power series expansion \eqref{appeq1} for $\sigma$ into $Q_1(s)$ and expanding the Bessel function, we obtain
\begin{align*}
   &Q_1(s,t)=\sum_{n=0}^{\infty} 2 p_n \sqrt{C (s-a) (t-c)} (s-a)^n \int_0^{\frac{\pi}{2}} \cos ^{2 n+1} \theta \sum_{m=0}^{\infty} \frac{[C(s-a) (t-c)]^{m+\frac{1}{2}}}{m!(m+1)!}\sin ^{2 m+1} \theta d \theta
\end{align*}
The interchange of the summation over $n$ and the integral is justified by the uniform convergence of the partial sums of \eqref{appeq1} on the interval $[a, b]$. Similarly, we exchange the order of the summation over $m$ and the integral, yielding:
\begin{align*}
    Q_1(s,t)&= \sum_{n=0}^{\infty} \sum_{m=0}^{\infty} 2 p_n \frac{[C(t-c)]^{m+1} (s-a)^{n+m+1}}{m!(m+1)!} \int_0^{\pi/2} \cos ^{2 n+1} \theta \sin ^{2 m+1} \theta d \theta \\
    & =\sum_{n=0}^{\infty} \sum_{m=0}^{\infty} p_n \frac{[C(t-c)]^{m+1} (s-a)^{n+m+1}}{(m+1)!} \frac{n!}{(n+m+1)!}\\
    & =\sum_{n=0}^{\infty} \sum_{m=1}^{\infty} p_n \frac{[C(t-c)]^{m} (s-a)^{n+m}}{m!} \frac{n!}{(n+m)!}
\end{align*}
where in the second line we applied the following identity
\begin{equation*}
    \int_0^{\pi / 2} \cos ^{2 n+1} \theta \sin ^{2 m+1} \theta=\frac{1}{2}\operatorname{Beta}(n+1, m+1)=\frac{\Gamma(n+1) \Gamma(m+1)}{2\Gamma(n+m+2)}
\end{equation*}
Similarly for $Q_2$ we obtain:
\begin{align*}
    Q_2(s,t) =\sum_{n=0}^{\infty} \sum_{m=1}^{\infty} q_n \frac{[C (s-a)]^{m} (t-c)^{n+m}}{m!} \frac{n!}{(n+m)!}
\end{align*}
Combining these results and expanding the remaining terms in Eq. \eqref{step1eq1}, we have
\begin{align*}
    k(s,t)&=-\sum_{m=0}^\infty p_0\frac{[C(t-c)]^m(s-a)^m}{(m!)^2} + \sum_{n=0}^\infty p_n(s-a)^n +\sum_{n=0}^{\infty} \sum_{m=1}^{\infty} p_n \frac{[C(t-c)]^{m} (s-a)^{n+m}}{m!} \frac{n!}{(n+m)!} \nonumber \\
    &+\sum_{n=0}^\infty q_n(t-c)^n+\sum_{n=0}^{\infty} \sum_{m=1}^{\infty} q_n \frac{[C (s-a)]^{m} (t-c)^{n+m}}{m!} \frac{n!}{(n+m)!} \nonumber\\
    &= \sum_{n=1}^{\infty} \sum_{m=0}^{\infty} p_n \frac{[C(t-c)]^{m} (s-a)^{n+m}}{m!} \frac{n!}{(n+m)!} + \sum_{n=0}^{\infty} \sum_{m=0}^{\infty} q_n \frac{[C (s-a)]^{m} (t-c)^{n+m}}{m!} \frac{n!}{(n+m)!}
\end{align*}

Given that $k(a,c)=\sigma(a)=\tau(c)$, which implies $p_0=q_0$,  Eq. \eqref{appeq2} can also be written as 
\begin{equation*}
    k(s,t)=\sum_{n=0}^{\infty} \sum_{m=0}^{\infty} p_n \frac{[C(t-c)]^{m} (s-a)^{n+m}}{m!} \frac{n!}{(n+m)!} + \sum_{n=1}^{\infty} \sum_{m=0}^{\infty} q_n \frac{[C (s-a)]^{m} (t-c)^{n+m}}{m!} \frac{n!}{(n+m)!}
\end{equation*}\qed

Proposition \ref{prop:a1} establishes the foundational framework necessary to prove both Proposition \ref{approxsolution} (also included here in the appendix as Proposition \ref{prop:a2}) and Proposition \ref{propinterp} (included here as Proposition \ref{prop:a3}). 
\begin{proposition}\label{prop:a2}
    Consider the Goursat PDE on the rectangle $[a,b]\times[c,d]$ as in \eqref{goursatonrectangle}, with initial boundaries condition given by \eqref{initialps}, then $\forall n\in\NN$ the coefficients of \eqref{finalps} are given by
    \begin{equation}\label{appeq111}
        \ptil_{n} = \sum_{k=0}^{n} p_k \frac{(C \dcd)^{n-k}\ k!}{(n-k)!\ n!} + \sum_{k=1}^{\infty} q_k \frac{ C^n\dcd^{n+k}\ k!}{(n+k)!\ n!}
    \end{equation}
    and, symmetrically
    \begin{equation}\label{appeq222}
        \qtil_{n} = \sum_{k=0}^{n} q_k \frac{(C \dab)^{n-k}\ k!}{(n-k)!\ n!} + \sum_{k=1}^{\infty} p_k \frac{ C^n\dab^{n+k}\ k!}{(n+k)!\ n!}
    \end{equation}
    where $\dab=(b-a)$ and $\dcd=(d-c)$.
\end{proposition}
\paragraph{Proof.} The proof follows directly. To establish \eqref{appeq111}, we set  $t = d$  in \eqref{appeq3} and then rearrange the sums in terms of  $(s-a)$. Similarly, by substituting  $s = c$  into \eqref{appeq2}, we can prove \eqref{appeq222}. \qed

\begin{proposition}\label{prop:a3}
Consider the Goursat problem on a rectangle $[a,b]\times[c,d]$ as in \eqref{goursatonrectangle}. Further, assume that the initial boundary conditions $k(\cdot,c):[a,b]\to\R$ and $k(a,\cdot):[c,d]\to\R$ are polynomials at $s=a$ and $t=c$ respectively, such that
\begin{equation}\label{eq:interpassumpproof}
    k(s,c)=\sum_{n=0}^N g_n (s-a)^n  \qquad k(a,t) =\sum_{n=0}^N h_n (t-c)^n
\end{equation}with the constraint $g_0=h_0$. Then the solution $k(s,t)$ can be written as
\begin{equation}\label{eqapp:interp}
    k(s,t) = g_0I_0(2\sqrt{C(s-a)(t-c)}) + \sum_{n=1}^{N} \left(g_n (s-a)^n + h_n(t-c)^n \right)\hypf(n+1;\, C(s-a)(t-c)) 
\end{equation}
with $s\in[a,b]$, $t\in[c,d]$. Here, $\hypf$ is the confluent hypergeometric limit function \cite{andrews1999special} defined by:
\begin{equation}\label{hyp0f1def appendix}
    \hypf(k;z)=\sum_{n=0}^\infty\frac{z^n}{(k)_nn!}
\end{equation} where $(k)_n$ is the Pochhammer symbol, which represents the rising factorial.
\end{proposition}
\paragraph{Proof.} Again the proof is straightforward. According to assumption \eqref{eq:interpassumpproof}, we focus on the truncation at level $n = N$ of the power series in \eqref{appeq3}, which rearranged yields:
\begin{align*}
    k(s,t)&= \sum_{n=0}^{N}p_n (s-a)^n \sum_{m=0}^{\infty}  \frac{[C(t-c)]^{m} (s-a)^{m}}{m!} \frac{n!}{(n+m)!} + \sum_{n=1}^{N} q_n(t-c)^n\sum_{m=0}^{\infty}\frac{[C (s-a)]^{m} (t-c)^{m}}{m!} \frac{n!}{(n+m)!}
\end{align*} 
Using Pochhammer notation, this simplifies to:
\begin{align}\label{appeq32}
    k(s,t)&= \sum_{n=0}^{N}p_n (s-a)^n \sum_{m=0}^{\infty}  \frac{[C(t-c)]^{m} (s-a)^{m}}{m!(n+1)_m}  + \sum_{n=1}^{N} q_n(t-c)^n\sum_{m=0}^{\infty}\frac{[C (s-a)]^{m} (t-c)^{m}}{m!(n+1)_m} 
\end{align} 
From the definition of the confluent hypergeometric function $\hypf$ - Eq. \eqref{hyp0f1def appendix}, we recover equation \eqref{eqapp:interp}. \qed

\newpage
\section{Additional Experiments}\label{addexp}

\begin{figure}[!h]
    \centering
    \includegraphics[width=0.9\textwidth]{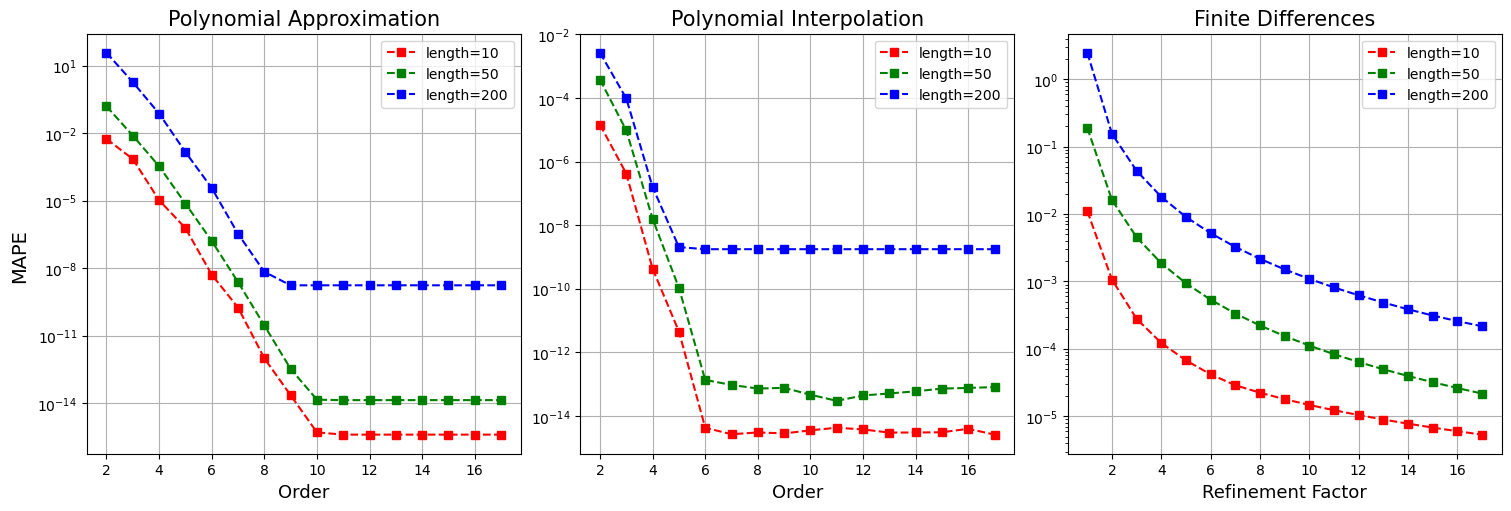}
    \caption{Comparison of MAPE against polynomial order for the polynomial approximation and interpolation schemes is shown on the \textbf{left} and \textbf{center}, respectively. On the \textbf{right}, MAPE is plotted against the refinement factor for the finite differences scheme. This analysis is performed on a $8\times8$ batch of 2-$d$ paths sampled from \eqref{sinpath} with different lengths (color-coded). }
    \label{fig:errorvslevelsin}
    \vspace{-0.1in}
\end{figure}

In this section, we present supplementary results to further evaluate the performance of our proposed numerical schemes, and, at the end, we provide tables reporting the numerical results corresponding to the experiment shown in Figure \ref{fig:errorvslevel}.  We begin by comparing the mean absolute percentage error (MAPE) of the schemes in a specific setting where we compute the schemes for a setting where we compute $k(x^\ell,y^\ell)$, with $x^\ell$ and $y^\ell$ being discrete stochastic processes defined as:
\begin{align}\label{sinpath}
    x_i = \sin(U_i)\qquad\textrm{and}\qquad y_i = \cos(V_i)
\end{align} where each $U_i,V_i$ are 2-$d$ vectors, and each entry is drawn independently from the standard normal distribution $\mathcal{N}(0,1)$ and $i\in\{1,\dots,\ell\}$.  The results, illustrated in Figure \ref{fig:errorvslevelsin}, demonstrate that our proposed methods consistently achieve significantly lower MAPE compared to the finite differences scheme introduced by \cite{salvi2021signature}. As in previous experiments, we use the truncated signature at level $21$ as a benchmark for accuracy.

\begin{figure}[!h]
    \centering
    \includegraphics[width=0.9\textwidth]{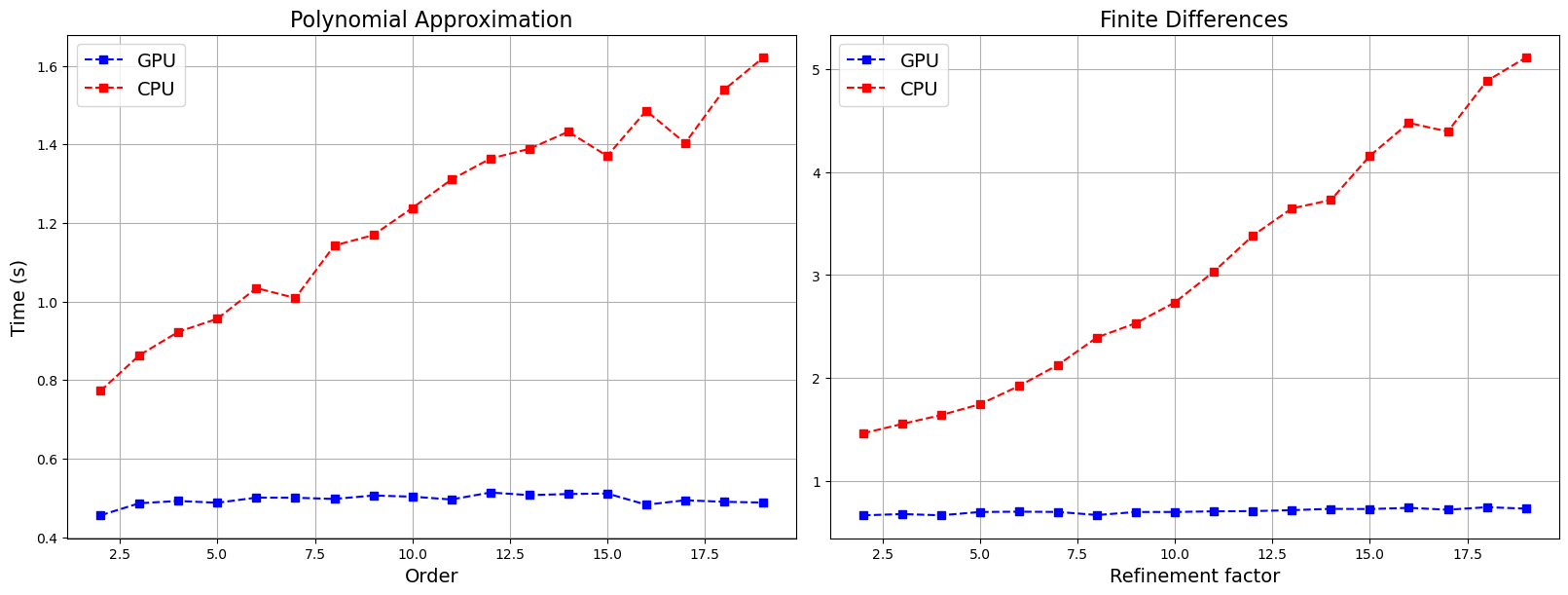}
    \caption{Comparison of elapsed time (in seconds) to compute the signature kernel between two batches of eight 8-$d$ Brownian motions, each of length 32, against the order approximation for the polynomial approximation scheme (\textbf{left}) and refinement factor for the finite differences scheme (\textbf{right}).}
    \label{fig:timeanal1}
\end{figure}

Figure \ref{fig:timeanal1} shows that for smaller datasets the polynomial approximation scheme exhibits near-linear complexity in $N$. This can be attributed to the vectorization capabilities of the CPU and the use of Single Instruction-Multiple Data (SIMD) instructions, which allow multiple coefficients to be processed simultaneously, reducing the effective computational load.

\begin{figure}[!h]
    \centering
    \includegraphics[width=0.9\textwidth]{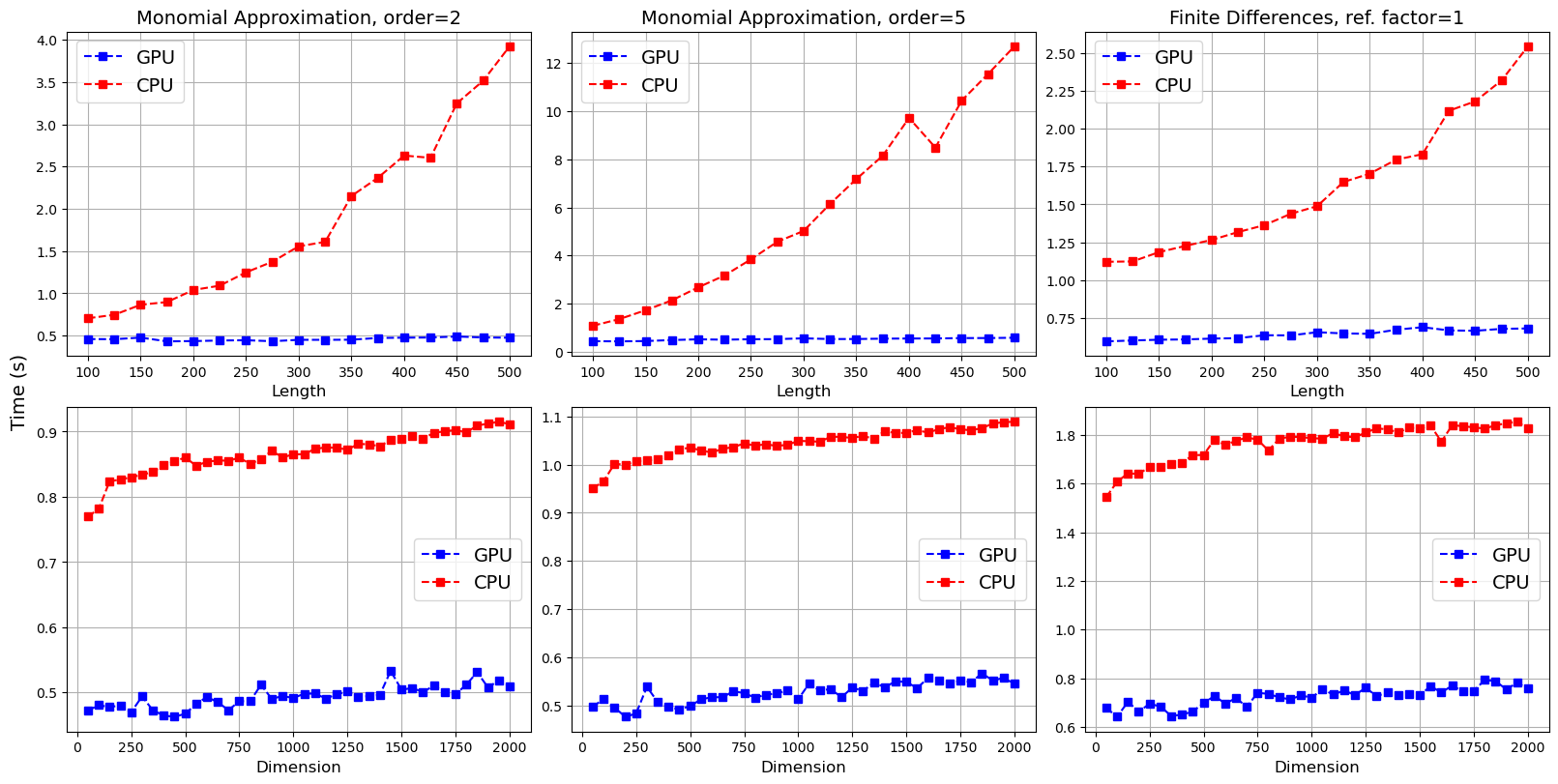}
    \caption{Comparison of elapsed time (in seconds) to compute the signature kernel between two batches of eight Brownian motions against length (\textbf{top row}) and dimension (\textbf{bottom row}).}
    \label{fig:timeanal2}
\end{figure}

In the context of GPU acceleration, Figure \ref{fig:timeanal2} (top row) illustrates how enabling parallelization breaks the quadratic complexity in relation to the time series length, reducing it to linear complexity. The bottom row highlights the computational impact of the time series dimension. A key distinction between the untruncated and truncated signature kernels lies in this aspect: while the truncated version requires computing the signature, whose size grows exponentially with the number of dimensions, solving the PDE does not face this limitation.

\newpage
$ $
\newpage

\begin{table}[]
    \centering
    \setlength{\tabcolsep}{12pt}  
    \renewcommand{\arraystretch}{1.1}  
    \caption{Numerical results for the mean absolute percentage error (MAPE) corresponding to the experiment illustrated in Figure \ref{fig:errorvslevel} (left) for paths of length 10.}
    \begin{tabular}{c c c c}
        \toprule
        Order/Ref. Factor & Poly. Approx. (\%) & Poly Interp. (\%) & Finite Differences (\%) \\
        \midrule
        1  & --                & --                & $1.103 \times 10^{-2}$ \\
        2  & $5.588 \times 10^{-3}$  & $1.456 \times 10^{-5}$  & $1.064 \times 10^{-3}$ \\
        3  & $7.188 \times 10^{-4}$  & $4.195 \times 10^{-7}$  & $2.785 \times 10^{-4}$ \\
        4  & $1.060 \times 10^{-5}$  & $4.283 \times 10^{-10}$ & $1.223 \times 10^{-4}$ \\
        5  & $6.135 \times 10^{-7}$  & $4.562 \times 10^{-12}$ & $6.745 \times 10^{-5}$ \\
        6  & $5.253 \times 10^{-9}$  & $4.291 \times 10^{-15}$ & $4.181 \times 10^{-5}$ \\
        7  & $1.812 \times 10^{-10}$ & $2.689 \times 10^{-15}$ & $2.888 \times 10^{-5}$ \\
        8  & $1.021 \times 10^{-12}$ & $3.056 \times 10^{-15}$ & $2.238 \times 10^{-5}$ \\
        9  & $2.303 \times 10^{-14}$ & $2.856 \times 10^{-15}$ & $1.795 \times 10^{-5}$ \\
        10 & $5.185 \times 10^{-16}$ & $3.557 \times 10^{-15}$ & $1.473 \times 10^{-5}$ \\
        11 & $4.118 \times 10^{-16}$ & $4.309 \times 10^{-15}$ & $1.231 \times 10^{-5}$ \\
        12 & $4.118 \times 10^{-16}$ & $3.798 \times 10^{-15}$ & $1.043 \times 10^{-5}$ \\
        13 & $4.118 \times 10^{-16}$ & $3.023 \times 10^{-15}$ & $8.954 \times 10^{-6}$ \\
        14 & $4.118 \times 10^{-16}$ & $3.047 \times 10^{-15}$ & $7.769 \times 10^{-6}$ \\
        15 & $4.118 \times 10^{-16}$ & $3.097 \times 10^{-15}$ & $6.804 \times 10^{-6}$ \\
        16 & $4.118 \times 10^{-16}$ & $3.979 \times 10^{-15}$ & $6.007 \times 10^{-6}$ \\
        17 & $4.118 \times 10^{-16}$ & $2.610 \times 10^{-15}$ & $5.343 \times 10^{-6}$ \\
        \bottomrule
    \end{tabular}
    \label{tab1}
\end{table}

\begin{table}[]
    \centering
    \setlength{\tabcolsep}{12pt}  
    \renewcommand{\arraystretch}{1.1}  
    \caption{Numerical results for the mean absolute percentage error (MAPE) corresponding to the experiment illustrated in Figure \ref{fig:errorvslevel} (center) for paths of length 50.}
    \begin{tabular}{c c c c}
        \toprule
        Order/Ref. Factor & Poly. Approx. (\%) & Poly Interp. (\%) & Finite Differences (\%) \\
        \midrule
        1  & --                & --                & $1.898 \times 10^{-1}$ \\
        2  & $1.567 \times 10^{-1}$  & $3.797 \times 10^{-4}$  & $1.611 \times 10^{-2}$ \\
        3  & $7.884 \times 10^{-3}$  & $9.637 \times 10^{-6}$  & $4.513 \times 10^{-3}$ \\
        4  & $3.415 \times 10^{-4}$  & $1.586 \times 10^{-8}$  & $1.863 \times 10^{-3}$ \\
        5  & $7.536 \times 10^{-6}$  & $1.072 \times 10^{-10}$ & $9.386 \times 10^{-4}$ \\
        6  & $1.641 \times 10^{-7}$  & $1.370 \times 10^{-13}$ & $5.358 \times 10^{-4}$ \\
        7  & $2.461 \times 10^{-9}$  & $9.547 \times 10^{-14}$ & $3.345 \times 10^{-4}$ \\
        8  & $2.903 \times 10^{-11}$ & $7.337 \times 10^{-14}$ & $2.222 \times 10^{-4}$ \\
        9  & $3.423 \times 10^{-13}$ & $7.727 \times 10^{-14}$ & $1.548 \times 10^{-4}$ \\
        10 & $1.472 \times 10^{-14}$ & $4.703 \times 10^{-14}$ & $1.120 \times 10^{-4}$ \\
        11 & $1.399 \times 10^{-14}$ & $3.010 \times 10^{-14}$ & $8.344 \times 10^{-5}$ \\
        12 & $1.403 \times 10^{-14}$ & $4.448 \times 10^{-14}$ & $6.375 \times 10^{-5}$ \\
        13 & $1.403 \times 10^{-14}$ & $5.139 \times 10^{-14}$ & $4.982 \times 10^{-5}$ \\
        14 & $1.403 \times 10^{-14}$ & $6.015 \times 10^{-14}$ & $3.965 \times 10^{-5}$ \\
        15 & $1.403 \times 10^{-14}$ & $7.312 \times 10^{-14}$ & $3.204 \times 10^{-5}$ \\
        16 & $1.403 \times 10^{-14}$ & $7.684 \times 10^{-14}$ & $2.624 \times 10^{-5}$ \\
        17 & $1.403 \times 10^{-14}$ & $8.157 \times 10^{-14}$ & $2.175 \times 10^{-5}$ \\
        \bottomrule
    \end{tabular}
    \label{tab2}
\end{table}

\begin{table}[]
    \centering
    \setlength{\tabcolsep}{12pt}  
    \renewcommand{\arraystretch}{1.1}  
    \caption{Numerical results for the mean absolute percentage error (MAPE) corresponding to the experiment illustrated in Figure \ref{fig:errorvslevel} (right) for paths of length 200.}
    \begin{tabular}{c c c c}
        \toprule
        Order/Ref. Factor & Poly. Approx. (\%)    & Poly Interp. (\%)  & Finite Differences (\%) \\
        \midrule
        1  & --              & --                 & $2.468 \times 10^0$ \\
        2  & $3.766 \times 10^{1}$   & $2.661 \times 10^{-3}$  & $1.546 \times 10^{-1}$ \\
        3  & $1.867 \times 10^{0}$   & $9.724 \times 10^{-5}$  & $4.366 \times 10^{-2}$ \\
        4  & $7.480 \times 10^{-2}$  & $1.588 \times 10^{-7}$  & $1.797 \times 10^{-2}$ \\
        5  & $1.000 \times 10^{-3}$  & $2.058 \times 10^{-9}$  & $9.044 \times 10^{-3}$ \\
        6  & $3.708 \times 10^{-5}$  & $1.769 \times 10^{-9}$  & $5.200 \times 10^{-3}$ \\
        7  & $3.418 \times 10^{-7}$  & $1.768 \times 10^{-9}$  & $3.252 \times 10^{-3}$ \\
        8  & $6.997 \times 10^{-9}$  & $1.770 \times 10^{-9}$  & $2.164 \times 10^{-3}$ \\
        9  & $1.762 \times 10^{-9}$  & $1.768 \times 10^{-9}$  & $1.509 \times 10^{-3}$ \\
        10 & $1.766 \times 10^{-9}$  & $1.768 \times 10^{-9}$  & $1.093 \times 10^{-3}$ \\
        11 & $1.767 \times 10^{-9}$  & $1.768 \times 10^{-9}$  & $8.151 \times 10^{-4}$ \\
        12 & $1.766 \times 10^{-9}$  & $1.769 \times 10^{-9}$  & $6.234 \times 10^{-4}$ \\
        13 & $1.766 \times 10^{-9}$  & $1.769 \times 10^{-9}$  & $4.868 \times 10^{-4}$ \\
        14 & $1.767 \times 10^{-9}$  & $1.770 \times 10^{-9}$  & $3.870 \times 10^{-4}$ \\
        15 & $1.766 \times 10^{-9}$  & $1.769 \times 10^{-9}$  & $3.124 \times 10^{-4}$ \\
        16 & $1.767 \times 10^{-9}$  & $1.770 \times 10^{-9}$  & $2.582 \times 10^{-4}$ \\
        17 & $1.766 \times 10^{-9}$  & $1.769 \times 10^{-9}$  & $2.167 \times 10^{-4}$ \\
        \bottomrule
    \end{tabular}
    \label{tab3}
\end{table}

\end{document}